\documentclass[12pt]{amsart}

\usepackage{latexsym}
\usepackage{amssymb}
\usepackage[centertags]{amsmath}

\usepackage{verbatim}

\newtheorem{theorem}{Theorem}
\newtheorem{lemma}[theorem]{Lemma}

\newtheorem{claim}[theorem]{Claim}
\newtheorem{corollary}[theorem]{Corollary}
\newtheorem{definition}[theorem]{Definition}

\numberwithin{theorem}{section}

\newcommand{\wh}{\widehat}

\renewcommand{\ge}{\geqslant}
\renewcommand{\le}{\leqslant}

\renewcommand{\phi}{\varphi}

\renewcommand{\#}{\operatorname{Card}}

\newcommand{\clos}{\operatorname{clos}}
\newcommand{\card}{\operatorname{card}}
\newcommand{\supp}{\operatorname{supp}}

\newcommand{\dist}{\operatorname{dist}}

\renewcommand{\Im}{\operatorname{Im}}

\newcommand{\la}{\lambda}
\renewcommand{\d}{{\rm d}}

\renewcommand{\epsilon}{\varepsilon}

\renewcommand{\le}{\leqslant}
\renewcommand{\ge}{\geqslant}
\renewcommand{\epsilon}{\varepsilon}

\newcommand{\E}{\mathcal E}

\newcommand{\R}{\mathbb R}
\newcommand{\C}{\mathbb C}
\newcommand{\Z}{\mathbb Z}

\title[Weighted exponential approximation]
{Weighted exponential approximation and non-classical
orthogonal spectral measures}

\thanks{
This research is a part of the ESF Networking Programme
"Harmonic and Complex Analysis and its Applications".
The first named author was partially supported by the
ANR projects DYNOP and FRAB}



\author{Alexander Borichev}

\address{A.B.: Centre de Math\'{e}matiques et Informatique\\
Universit\'{e} Aix-Marseille \\
39 rue Joliot-Curie\\
13453 Marseille, Cedex 13\\
France}

\email{borichev@cmi.univ-mrs.fr}

\author{Mikhail Sodin}

\address{M.S.: School of Mathematics\\
Tel Aviv University\\
Tel Aviv 69978\\
Israel}

\email{sodin@post.tau.ac.il}

\begin{document}

\bigskip
\dedicatory{\hfill To Vladimir Marchenko with admiration}

\begin{abstract}
A long-standing open problem in harmonic analysis is: {\em given
a non-negative measure $\mu$ on $\R$, find the infimal width of frequencies
needed to approximate any function in $L^2(\mu)$}. We consider
this problem in the ``perturbative regime'', and characterize asymptotic
smallness of perturbations of measures which do not change that infimal
width. Then we apply this result to show that there are no local
restrictions on the structure of orthogonal spectral measures of
one-dimensional Schr\"odinger operators on a finite interval.
This answers a question raised by V.~A.~Marchenko.
\end{abstract}

\maketitle

\section{Introduction and main results}\label{sect_Intro}

\subsection{The type problem}\label{subsect_1.1}
We say that a non-negative measure $\mu$ on $\R$ has at most polynomial
growth if, for some $s<\infty$,
\begin{equation}\label{eq.growth1}
\int_\R \frac{\d\mu(\la)}{1+|\la|^{2s}} <\infty\,.
\end{equation}
For a measure $\mu$ of at most polynomial growth, we define its {\em type}
$T(\mu)$ as follows. Let $\E(a)=\mathfrak F C_0^\infty(-a,a)$ be
the Fourier image of the space of the $C^\infty$-smooth complex-valued functions compactly
supported by $(-a,a)$. This is the linear space of entire functions of
exponential type less than $a$ that decay on the real axis faster than any
negative power of $|\la|$. Then $\E(a)\subset L^2(\mu)$, and
\begin{multline*}
T(\mu) \stackrel{\rm def}= \inf\,\bigl\{ a\colon \E(a) {\rm\ is\ dense\ in
\ } L^2(\mu) \bigr\}\\
= \sup\, \bigl\{ a\colon \E(a) \text{\rm\ isn't dense in\ } L^2(\mu) \bigr\}
\,,
\end{multline*}
which  is one half of the infimal width of the spectrum
needed to approximate any function in $L^2(\mu)$.

The definition of the type is not sensitive to the choice of a ``natural''
linear space $\E(a)$ of entire functions of exponential type at most $a$.
For instance, it is not difficult to check that if the Paley--Wiener space
$\mathfrak F L^2(-a,a)$ is contained in $L^2(\mu)$, then without affecting
the definition of the type, one can replace therein the linear space $\E(a)$
by $\mathfrak F L^2(-a,a)$. If the measure $\mu$ is finite, then one can replace $\E(a)$
by the space of the finite linear combinations of exponential functions
$\la\mapsto e^{{\rm i}t\la}$ with $-a< t< a$.

The range of $T(\mu)$ is $[0,\infty]$ with both ends included. If the tails
of the measure $\mu$ decay so fast that the polynomials belong to the space
$L^2(\mu)$ and are dense therein, then it is easy to see that $T(\mu)=0$.
Another instance of the zero type occurs when the support of the measure $\mu$
has long gaps. On the other hand, Lebesgue measure $m$ on the real axis has
infinite type. An intermediate case occurs for the sum of point masses
$\delta$ at arithmetic progression: for $0<\ell<\infty$, we have
$T\bigl(\sum_{\la\in \ell\mathbb Z}\delta_\la\,\bigr)=\pi\ell^{-1}$.
These are toy models for other situations when the type $T(\mu)$ can be
explicitly computed. For reader's orientation, we bring a short summary of
what is known in Appendix~\ref{appendix-1}.

For any $g\in L^2(\mu)$, the Fourier transform of the measure $g\,\d\mu$
is a tempered distribution. If $\E(a)$ is not dense in $L^2(\mu)$ and
$g\in L^2(\mu)$ is orthogonal to $\E(a)$, then the Fourier transform of
$g\d\mu$ vanishes on $(-a,a)$. In this case, we say that the interval
$(-a,a)$ is a spectral gap of $g$. Therefore, {\em the type $T(\mu)$
coincides with one half of the supremal length of spectral gaps of functions
in $L^2(\mu)$}.

\subsection{Relations with other classical problems of analysis}

The problem of {\em effective computation of the type $T(\mu)$}, for brevity,
the {\em type problem}, is intimately related to other classical problems in
analysis. It originates from the works of Kolmogorov and Wiener on the
prediction of Gaussian stationary processes, see \cite{Krein45},
\cite[\S~3.7 and Chapter~4]{DM}. Then Gelfand and Levitan~\cite[\S8]{GL} and
Krein~\cite[Theorem~4]{Krein53} and~\cite{Krein54} discovered a deep
relation between the type problem and the spectral theory of the second
order ordinary differential operators which we discuss in
Section~\ref{subsection-inverse_spectral}.

The type problem is one of the central problems in the de~Branges theory of Hilbert
spaces of entire functions. By this theory, given a non-negative measure
$\mu$ satisfying the property
$$
\int_{\R} \frac{\d\mu (\la)}{1+\la^2} < \infty\,,
$$
there exists a unique chain of de Branges Hilbert
spaces of entire functions $\mathcal H(E_t)$ such that

\smallskip\par\noindent (i) the entire functions $E_t$ are of Cartwright class\footnote{
That is, the functions $E_t$ have at most exponential type and convergent
logarithmic integral $\int_\R \frac{\log^+|E_t(x)|}{1+x^2}\, {\rm d} x < \infty $.
},

\smallskip\par\noindent
(ii) $\mathcal H(E_{t_1})$ is
contained
isometrically in $\mathcal H(E_{t_2})$ for $t_1<t_2$,

\smallskip\par\noindent
(iii) each space
$\mathcal H(E_t)$ is contained isometrically in $L^2(\mu)$,

and

\smallskip\par\noindent
(iv)
$\bigcup_t \mathcal H(E_t)$ is dense in $L^2(\mu)$.

\par\noindent Then it is not difficult
to show that the supremum of exponential types of the functions $E_t$
coincides with $T(\mu)$. There is a remarkable formula due to
Krein~\cite{Krein54} and de~Branges~\cite[Theorem~39]{deBr} that expresses
the type $T(\mu)$ via the coefficients of the second order canonical system
describing evolution along the chain
of the spaces $\mathcal H(E_t)$. Note that de~Branges' book~\cite{deBr}
contains a wealth of results (Theorems~61--68) related to the type problem.

It is also worth mentioning that the type problem is a part of the general
Bernstein weighted approximation problem, cf. Dym~\cite{Dym},
Pitt~\cite{Pitt}, Koosis~\cite[Chapters~VI-VII]{Koosis1}, and
Levin~\cite{Levin}. The methods developed for solving the Bernstein problem will be
used in this work.

At last, the type problem is closely connected with fundamental results
of Beurling and Malliavin on multipliers and the radius of completeness,
see Koosis \cite{Koosis2}, and recent works of Mitkovski and
Polto\-ratski~\cite{PolMit} and of Poltoratski~\cite{Pol}. The papers by
Mitkovski and Polto\-ratski suggest a novel approach to the type problem based
on injectivity of certain Toeplitz operators.

\subsection{Perturbations of measures}

Since we do not know how to compute the type, it is natural to ask
which perturbations of positive measures {\em preserve} their types?
We prove that exponentially small perturbations of measures do not change
their types and then we show that this result is sharp.

Given $\delta>0$ and $x\in\R$, we denote
\begin{equation}
I_x=I_{x,\delta}=[x-e^{-\delta|x|}, x+e^{-\delta|x|}].
\label{m0}
\end{equation}
By $k I_x=k I_{x,\delta}$ we denote the closed interval centered at $x$ with length $k$
times that of $I_x$.

\begin{definition}[majorization in mean with exponentially small
error]\label{def_maj}
We write $\mu \preccurlyeq \widetilde\mu$ if there exist constants
$\delta>0$, $C>0$, and $n\ge 0$, such that, for all $x\in\R$,
\[
\mu(I_{x,\delta})\le C(1+|x|)^n\left(\widetilde\mu(2I_{x,\delta})+
e^{-2\delta|x|} \right).
\]
\end{definition}

\begin{definition}[stable density]\label{def-stable_density}
We say that $\E(a)$ is stably dense in $L^2(\mu)$ if for each $t\ge 0$,
$\E(a)$ is dense in $L^2(\mu_t)$ with
${\rm d}\mu_t (\la) = (1+|\la|)^t\, {\rm d}\mu(\la)$.
\end{definition}
It is not difficult to show (see Appendix~\ref{appendix-2})
that $\E(a)$ is stably dense in $L^2(\mu)$
if and only if, for any finite set of points $\la_1, ..., \la_N\in\R$,
$\E(a)$ is dense in $L^2(\widetilde{\mu})$, where $\widetilde{\mu} = \mu
+ \sum_k \delta_{\la_k}$. Here and everywhere below, $\delta_\la$ is a
unit point mass at $\la$.

Non-stable density is quite non-generic, though it often occurs in various
applications. It is known that if $\E(a)$ is dense but
not stably dense, then the measure $\mu$ is supported by the zero set of an
entire function of exponential type $a$, see Lemma~\ref{lemma-singular} in
Appendix~\ref{appendix-2} for a more precise statement.

\begin{theorem}\label{thm_main1}
Let $\widetilde\mu$ be a non-negative measure of at most polynomial growth,
and let $\mu\preccurlyeq\widetilde\mu$. If $\E(a)$ is stably dense in
$L^2(\widetilde{\mu})$, then $\E(a)$ is dense in $L^2(\mu)$.
\end{theorem}

By the remark made before the theorem, if $\E(a)$ is dense in $L^2(\mu)$,
then for each $a'>a$, $\E(a')$ is stably dense in $L^2(\mu)$. Hence,

\begin{corollary}\label{cor-main}
Let $\widetilde\mu$ be a non-negative measure of at most polynomial growth,
and let $\mu\preccurlyeq\widetilde\mu$. Then $T(\mu)\le T(\widetilde{\mu})$.
\end{corollary}

\medskip
We call the measures $\widetilde\mu$ and $\mu$ {\em weakly equivalent} if
$\widetilde\mu\preccurlyeq\mu$ and $\mu\preccurlyeq\widetilde\mu$.
Note that if two positive measures coincide outside of a finite interval,
then they are weakly equivalent.

\begin{corollary}\label{cor_Equiv}
Weakly equivalent measures have equal types.
\end{corollary}

We note that Theorem~\ref{thm_main1} has a counterpart, which deals with
polynomial density in $L^2(\mu)$, see Section~\ref{add1}.

\medskip
The following result shows that the statements of Corollaries~\ref{cor-main}
and~\ref{cor_Equiv} are sufficiently sharp:

\begin{theorem}\label{thm-sharpness}
Given a positive function $\varepsilon$ such that $\varepsilon(r)\to 0$,
$r\to\infty$,

\smallskip\par\noindent {\rm (i)} there exists a function $\varphi$ such
that $T\big(\varphi(x)\,\d x\big)=0$ while \newline
$T\big((\varphi(x)+e^{-\varepsilon(|x|)|x|})\,\d x\big)=\infty$.

\smallskip\par\noindent {\rm (ii)} there exist two sequences of points
$\{x_n\}_{n\in\Z}$, $\{y_n\}_{n\in\Z}$,
$|y_n-x_n|\le e^{-\varepsilon(|x_n|)|x_n|}$, such that
$T\big(\sum_{n\in\Z}\delta_{x_n}\big)=\pi$ while 
$T\big(\sum_{n\in\Z}(\delta_{x_n}+\delta_{y_n})\big)=2\pi$.
\end{theorem}

\medskip
We note that the situation changes if we consider perturbations of
sufficiently regular measures. The types of such measures are more stable,
see a corollary to a classical result of Duffin and Schaeffer cited in
Appendix~\ref{xx3} (perturbations of Lebesgue measure), and
Benedicks~\cite[Theorems 8]{Ben} (perturbations of the sum of point masses
at arithmetic progression).

\subsection{Spectral theory of one-dimensional Schr\"odinger operators}
\label{subsection-inverse_spectral}

First, we recall some classical facts pertaining to the spectral
theory of one-dimensional Schr\"odinger operators. Below, we follow the
first two chapters of~\cite{Marchenko} (see also~\cite[Chapter~1]{LS}
and~\cite{GL, Krein53}).

\subsubsection{A piece of Weyl's spectral theory}\label{subsubsect_Weyl}

Given $a$, $0<a\le \infty$, consider the Sturm--Liouville equation
\begin{equation}\label{eq_SL}
-y''+q(x)y=\la^2 y\,,\qquad 0\le x<a,
\end{equation}
with a real-valued potential $q\in C[0,a)$. Note that we use the
``momentum'' $\la$ (not the ``energy'' $\la^2$) as the spectral parameter,
and that we do not impose any restrictions on $q$ at the right end-point
$x=a$.

Take the solution $\omega (\la, x)$ satisfying the boundary condition
\begin{equation}\label{eq_BC}
y(0) = 1, \qquad y'(0) = h\,.
\end{equation}
For each $x\in[0,a)$, this is an entire function of $\la$. It satisfies the
integral Sturm--Liouville equation
\[
\omega(\la,x)=\cos\la x+h\,\frac{\sin \la x}{\la}
+ \int_0^x \frac{\sin \la (x-t)}{\la} q(t) \omega (\la, t)\, \d t\,,
\]
which easily yields the estimate
$$
\bigl| \omega (\la, x) - \cos \la x\bigr| \le e^{x|{\rm Im}\la|}\, \frac
{Q(x) + |h|}{|\la|-Q(x)}\,,
$$
with $Q(x)=\int_0^x|q|$ and $|\la|>Q(x)$. Thus, given $x$, the function
$\la\mapsto\omega(\la,x)$ is an entire function of exponential type $x$
bounded on the real axis.

Consider the Weyl transform
\[
\mathfrak W f(\la) = \int_0^a f(x) \omega (\la, x)\, \d x\,.
\]
This transform is well-defined on the subspace $L^2_0(0,a)$ of $L^2(0,a)$
consisting of the functions that vanish on a neighbourhood of the end point
$x=a$. Note that if $f=-u''+qu$, and
\begin{equation}
u'(0)=hu(0),
\label{sq9}
\end{equation}
then $\mathfrak Wf(\la)=\la^2\mathfrak Wu(\la)$. A celebrated theorem
of Weyl says that

\smallskip\par\noindent $\bullet$ {\em there
exists a measure $\mu$ supported by $\R\cup{\rm i}\R$ and symmetric with
respect to the origin, such that}
\[
\|f\|_{L^2(0,a)}=\|\mathfrak W f\|_{L^2(\mu)}\,,\qquad f\in L^2_0(0,a).
\]

\smallskip\noindent The measure $\mu$ is called a {\em spectral measure}
of the Sturm--Liouville problem~\eqref{eq_SL}--\eqref{eq_BC}.
The map $\mathfrak W$ extends to the isometry $L^2(0,a)\to L^2(\mu)$,
and the inverse map, defined by
\[
f(x) = \int_{\R \cup {\rm i}\R}
\mathfrak W f(\la)\, \omega (\la, x)\, \d\mu(\la)\,,
\]
is called the eigenfunction expansion associated with the Sturm--Liouville
problem~\eqref{eq_SL}--\eqref{eq_BC}.
If the image
$\mathfrak W L^2(0,a)$ spans the closed subspace $L^2_{\tt e}(\mu)$ of even
functions in $L^2(\mu)$, then the spectral measure $\mu$ is called
{\em orthogonal} (or sometimes, {\em principal}).
It is known that each Sturm--Liouville problem has orthogonal spectral measures.

For the reader's orientation, we mention that there is a one-to-one correspondence
between orthogonal spectral measures $\mu$ and self-adjoint extensions
to a dense subset of $L^2(0, a)$ of the
operator $-y''+q(x)y$ with boundary condition
\eqref{sq9}.
Each self-adjoint extension of this type
is unitarily equivalent to the operator of multiplication by $\la^2$ in $L^2_{\tt e}(\mu)$.
In the limit-point case at the end point $x=a$, when the self-adjoint extension is unique,
the operator has only one spectral measure, and it is orthogonal. In the limit-circle
case, there are many spectral measures and some of them are orthogonal,
while the others correspond to self-adjoint operators defined on
an extension of the space $L^2(0, a)$. See Akhiezer and Glazman~\cite[Appendix~II]{AG}
for the details of this correspondence.

Note that it follows from Weyl's theory that
\begin{equation}\label{eq:growth}
\int_{\R \cup {\rm i}\R} \frac{\d\mu(\la)}{1+|\la|^2} < \infty\,.
\end{equation}

\subsubsection{A piece of the theory developed by Gelfand--Levitan, Krein,
and Marchen\-ko}\label{subsection-GLKM}

Weyl proved his theorem in 1909--10. Forty years later, Gel\-fand--Levitan, Krein,
and Marchenko developed a beautiful theory that fully describes spectral
measures of the one-dimensional Schr\"odin\-ger operator, and tells how to recover
the potential $q$ from the spectral measure $\mu$; see Marchenko~\cite
{Marchenko2} for a very illuminating account of the development of this
theory.

Given a measure $\mu$ supported by $\R\cup{\rm i}\R$, symmetric with respect
to the origin and satisfying the growth condition~\eqref{eq:growth},
we define the function
\begin{equation}\label{eq_ConMom}
\Phi(x)=\Phi[\mu](x)=\int_{\R\cup{\rm i}\R}\frac{1-\cos\la x}{\la^2}\,
\d\mu(\la)\,.
\end{equation}
This transform was introduced and studied by Povzner and Krein in the
1940-s; for its basic properties see
\cite[items~10-12 in Addenda to Chapter~V]{Akh}\footnote{
Therein, the transform is written in the form
\[
\int_\R \frac{1-\cos (x\sqrt{s})}{s}\,\d\rho(s)\,,
\qquad {\rm with} \quad \d\rho(s)=2\d\mu(\sqrt{s})\,.
\]
}.

\begin{theorem}[Gelfand--Levitan]\label{thm-GL}
The measure $\mu$ is a spectral measure of the Sturm--Liouville boundary
problem~\eqref{eq_SL}--\eqref{eq_BC} with a continuous potential $q$ if and
only if
\begin{equation}\label{eq-central}
\Phi\in C^3 [0, 2a) \quad{\rm and} \quad
\Phi'(+0) = 1,\, \Phi''(+0)= -h\,.
\end{equation}
Moreover, the potential $q$ has the same number of continuous derivatives
on $[0,a)$ as $\Phi^{'''}$ has on $[0,2a)$.
\end{theorem}

\smallskip
Following~\cite[Chapter~2, \S 4]{Marchenko}, we rewrite the Gelfand--Levitan
condition~\eqref{eq-central} in a different form replacing the
$\Phi$-transform by the Fourier transform.
Let $\mu=\mu_\R+\mu_{{\rm i}\R}$, where the measure $\mu_\R$ is supported by
$\R$, and the measure $\mu_{{\rm i}\R}$ is supported by
${\rm i}\R\setminus\{0\}$. The r\^{o}les of these measures are very different.
The measure $\mu_\R$ is close in some sense to Lebesgue measure
$\tfrac1{\pi} m$, while the tails of the measure $\mu_{{\rm i}\R}$ decay
exponentially and $\mu_{{\rm i}\R}$ can be considered as a
``perturbation'' of the measure $\mu_\R$. More precisely, the
Gelfand--Levitan condition~\eqref{eq-central} is {\em equivalent} to
the following two conditions:

\smallskip\par\noindent {(GL--{\bf i})} there exists an even function
$M\in C^1(-2a,2a)$ such that the restriction of the distributional Fourier
transform\footnote{
We always use the following normalization
\[
\wh f (\la) = \int_\R f(x) e^{-{\rm i}\la x}\, \d x\,.
\]
for the Fourier transform.
}
$\widehat{\mu_\R}$ to $(-2a,2a)$ is equal to $2\delta_0+M$, that
is,
\begin{equation}
\int_\R\widehat{f}(\la)\,\d\mu(\la)=2f(0)+\int_{-2a}^{2a} f(x) M(x)\, \d x
\label{eq:add7}
\end{equation}
for each $f\in C^\infty_0(-2a, 2a)$;

\smallskip\par\noindent {(GL--{\bf ii})} the tails of the measure
$\mu_{{\rm i}\R}$ decay exponentially fast: for each $x<2a$,
\[
\int_0^\infty e^{x\la}\, \d\mu_{{\rm i}\R}({\rm i}\la) <\infty\,.
\]

\noindent If conditions (GL--{\bf i}) and (GL--{\bf ii}) are fulfilled, then
$h=-\widehat{\mu}(+0)$ (in other words,
$M(0)+\mu_{{\rm i}\R}({\rm i}\R)=-h$).

\smallskip
For instance, in the case of the zero potential on the semi axis ($q=0$,
$a=\infty$ and $h=0$), we have
\[
\mu_\R=\tfrac1{\pi}m,\qquad \mu_{{\rm i}\R}=0,\qquad \widehat{\mu}=2\delta_0,
\]
where $m$ is Lebesgue measure on $\R$, and $\Phi(x)=|x|$, $M(x)=0$, $x\in\R$.
In the case of the zero potential on a finite interval $[0,a)$ and $h=0$, we have
\[
\mu_\R=\tfrac 1{a}\sum_{n\in \Z} \delta_{\pi n/a},
\qquad \mu_{{\rm i} \R} = 0, \qquad
\widehat{\mu}=2\sum_{n\in\Z}\delta_{2n a},
\]
and $\Phi(x)=|x|$, $M(x)=0$ for $|x|<2a$.

\smallskip
The original statement of Theorem~\ref{thm-GL} in the paper by
Gelfand and Levitan~\cite{GL} contained a gap in one derivative between the necessary
and sufficient conditions which was removed by Krein in~\cite{Krein53}. Also
note that sometimes the Gelfand--Levitan theorem is formulated with an
additional assumption on the measure $\mu$ (e.g., see~\cite[Theorem~2.3.1]{Marchenko}):

\smallskip\par\noindent{(GL--{\bf iii})} for any $0<b<a$,
and for any even $f\in L^2(-b,b)$,
we have $\int |\widehat {f}|^2 \d \mu > 0$ unless $f=0$.

\smallskip\par\noindent However, Yavryan~\cite{Yavryan}
showed that this additional assumption (GL--{\bf iii})
follows\footnote{
Indeed, suppose that $\int |\widehat {f}|^2\d \mu = 0$.
Take any even continuous function $ \widetilde{M} $ with support in $(-2a,2a)$ coinciding with
$M$ in a neighborhood of $[-2b,2b]$, and denote by $\widetilde{\mu}$ a (signed) measure supported by $\R$ whose
distributional Fourier transform coincides with $2\delta_0 + \widetilde{M}$.
Then $d\widetilde{\mu}(t)/dt\to 2$, $|t|\to\infty$, and it is not difficult to see that for some $\varepsilon>0$ and for every
$g\in \mathfrak F L^2(-\varepsilon,\varepsilon)$ we have
$\int |\widehat {f}|^2|g|^2\,\d \widetilde\mu = 0$.
We apply this equality with $g(z)=h(\frac zA-1)$,
where $h(z)=\bigl((\sin z)/z\bigr)^{2N}$,
$A=2N/\varepsilon$, for large $N\in\mathbb N$. Then for some positive $c$ and for sufficiently large $A$, we have
$\int_{A-1}^{A+1}|\widehat {f}(t)|^2\,
\d \widetilde\mu(t)\le \exp(-cA)$.
Therefore, $|\widehat {f}(t)|\le c_1\exp(-ct)$, $t\to\infty$.
Since $\widehat {f}$ is in the Cartwright class, we get $f=0$.}
from condition (GL--{\bf i}).

\subsubsection{Orthogonal spectral measures}

Now, we turn to the orthogonality condition. Denote by
$\mathfrak CL^2_0(0,a)$ the image of the linear space $L^2_0(0,a)$ under the
cosine-transform
\[
\mathfrak C f(\la)=\int_0^a f(x)\cos\la x\,\d x\,,\qquad f\in L^2_0(0,a)\,.
\]
By the Paley--Wiener theorem, $\mathfrak CL^2_0(0,a)$ coincides with the
linear space of even entire functions of exponential type less than $a$
which belong to the space $L^2(\R)$.
The images of $ L^2_0(0,a) $ under the cosine-transform and
the Weyl transform coincide; i.e., $\mathfrak W L^2_0(0, a) =\mathfrak  CL^2_0(0,a) $
as linear spaces of entire functions. This follows from the classical
equations
\begin{align*}
\omega (\la, x)&=\cos \la x+\int_0^x K(x, t; h) \cos \la t {\rm d} t,\\
\cos \la x &=\omega(\la, x)+\int_0^x L(x, t; h) \omega(\la, t) {\rm d} t,
\end{align*}
with the kernels $ K $ and $ L $ continuous for
$ 0 \le x, t < a $, cf. \cite[Chapter~2, \S 2]{Marchenko}.
Since $L^2_0(0,a)$ is dense in $L^2(0,a)$ and since $\mathfrak W$
is an isometry between $L^2(0,a)$ and $L^2(\mu)$, we conclude that

\smallskip\par\noindent $\bullet$ {\em the spectral measure $\mu$ is
orthogonal if and only if $\mathfrak C L^2_0(0,a)$ is dense in the subspace
of even functions $L^2_{\tt e}(\mu)$}.

\smallskip
The latter condition follows from the
density of $\E(a)$ in the whole space $L^2(\mu)$. Hence, {\em the spectral
measure $\mu$ is orthogonal provided that $\E(a)$ is dense in $L^2(\mu)$}.
This relates orthogonal spectral measures to a weighted
exponential approximation problem, though a peculiar one: now, $\mu$ is a
symmetric measure supported by $\R \cup {\rm i}\R$.

\subsection{Non-classical orthogonal spectral measures}

If $a<\infty$ and the potential $q$ is continuous at the right end point
$x=a$, we arrive at the classical Sturm--Liouville problem with two regular
end points.
In this case, the measure $\mu$ is discrete and has a well-known
asymptotic behavior. For instance, if $a=\pi$, then
$\mu_\R=\sum_{n\in\Z}\alpha_n\delta_{\la_n}$ with $\la_n = n + O\bigl( \tfrac1n\bigr)$,
$\la_{-n}=\la_n$, and $\alpha_n = \tfrac1{\pi} + O\bigl( \tfrac1n \bigr)$,
$\alpha_{-n}=\alpha_n$, while $\mu_{{\rm i}\R}$ may consist only of finitely
many atoms. However, in the general case, when $x=a$ is a finite singular
end point, the situation is more tangled, and not much is known about
orthogonal spectral measures.
Even the most basic question: {\em what are
the restrictions imposed on the local structure of spectral measures of
Schr\"odinger operators on a finite interval by the orthogonality
condition?} remained open\footnote{
In the case $a=\infty$ there are {\em no} restrictions on the
local structure of orthogonal spectral measures, see
Gelfand and Levitan~\cite[\S~8]{GL}. In this respect, the case of the infinite
interval $[0, \infty)$ is much simpler since in that case {\em any} spectral measure
supported by $\R$ is automatically orthogonal. This follows from the density of
the space
$\bigcup_{a<\infty}\,\mathfrak C L^2_0(0, a)$ of even entire functions
of exponential type that belong to $L^2(\R)$ in the space of even functions
$L^2_{\tt e}(\mu)$ where $\mu$ is a measure
supported by $\R$ and satisfying \eqref{eq.growth1}.
}
since the 1950-s, see Marchenko~\cite{Marchenko2}.
The only result in this direction we are aware of is a delicate
construction by Pearson~\cite{Pearson}. He builds a potential $q$ on a finite
interval for which the orthogonal spectral measure $\mu$ is absolutely
continuous on a finite interval $[-\la_0,\la_0]$, the restriction of $\mu$
to $\R\setminus[-\la_0,\la_0]$ is discrete, and the restriction of $\mu$
to ${\rm i}\R\setminus\{0\}$ is at most the sum of finitely many atoms.
Pearson writes: ``The existence of an absolutely continuous spectrum for a
Schr\"odinger operator in a finite interval may be regarded as an exceptional
phenomenon, and we have to work quite hard to achieve it''~\cite[p.495]{Pearson-book}.
Curiously enough, Theorem~\ref{thm_main1} tells us that from the point of view of
the inverse spectral theory, this phenomenon is not exceptional.

\begin{definition}[stably orthogonal spectral measures]
\label{def-stably-orthogonal-measures}
We call a spectral measure $\mu$ of the Sturm--Liouville
problem~\eqref{eq_SL}--\eqref{eq_BC} on a finite interval
$[0,a)$ stably orthogonal if $\E(a)$ is stably dense in $L^2(\mu_\R)$.
\end{definition}
\noindent The orthogonal spectral measures corresponding to the classical
Sturm-Liouville problems with two regular end-points, or more generally, to
the problems with the limit-circle case at $x=a$, are not stable.
In Appendix~\ref{appendix-5} we bring an explicit
construction of a rather wide class of discrete spectral
measures, and in Appendix~\ref{appendix-6} we show that some of
them are stably orthogonal. We perturb
these measures, proving that there are
{\em no} local restrictions on the structure of orthogonal spectral measures
of the Sturm--Liouville problem on a finite interval:

\begin{theorem}\label{thm.main}
Suppose that $\mu_0$ is a stably orthogonal spectral measure of a
Sturm--Liouville problem~\eqref{eq_SL}--\eqref{eq_BC} on a finite interval
$[0,a)$, and that the measure $\mu_0$ is supported by the real axis. Suppose
that $\mu=\mu_\R+\mu_{{\rm i}\R}$ is a non-negative symmetric measure on
$\R\cup{\rm i}\R$ such that

\smallskip\par\noindent {\rm (i)} the support of the measure $\mu_\R$ does
not coincide with the zero set of an entire function of exponential
type $\le a$;

\smallskip\par\noindent {\rm (ii)} the integral
$\displaystyle\int_0^\infty \d (\mu_\R-\mu_0)$ converges, and
for some $\delta>0$ we have
\begin{equation}\label{eq_proximity}
\int_0^\infty e^{\delta\la } \left| \int_\la^\infty \d (\mu_\R - \mu_0) \right| \d \la < \infty\,,
\end{equation}
and {\rm (iii)}
\begin{equation}\label{eq.imag}
\int_0^\infty e^{\delta \la^2}\, \d\mu_{{\rm i}\R}({\rm i}\la) < \infty\,.
\end{equation}
Then the measure $\mu$ is an orthogonal spectral measure of a
Sturm--Liouville problem~\eqref{eq_SL}--\eqref{eq_BC} on the same interval
$[0,a)$, with the potential $q$ of the same class of smoothness on $[0,a)$
as that of the potential $q_0$ that corresponds to $\mu_0$.
\end{theorem}

Let us comment on the r\^{o}le of the technical condition (i) in Theorem~\ref{thm.main}.
Condition~\eqref{eq_proximity} allows us to apply Theorem~\ref{thm_main1} and to conclude
that $\E(a)$ is dense in $L^2(\mu_\R)$. Then condition (i) will allow us to
apply Lemma~\ref{lemma-singular}
from Appendix~\ref{appendix-2} and conclude that $\E(a)$ is stably
dense in $L^2(\mu_\R)$. We need this conclusion in order to further perturb
$\mu_\R$ by a measure $\mu_{i\R}$ supported by the imaginary axis.

\subsection{A brief reader's guide}
The rest of the paper consists of three parts, mostly
independent of each other, and four appendices. In the
first part (Sections~\ref{se2} and~\ref{se3}), we prove Theorem~\ref{thm_main1} which
says that exponentially
small perturbations of measures do not change their type.
In Section~\ref{se2}, we recall necessary preliminaries, and the proof
itself occupies Section~\ref{se3}. In the second part (Section~\ref{se4}),
we construct examples that show how sharp is Theorem~\ref{thm_main1}.
In the third part (Section~\ref{section-4}) we turn to perturbations
of stably orthogonal spectral measures and prove Theorem~\ref{thm.main}.
At the end of the paper we give information on the cases when the type can be explicitly calculated (Appendix A),
on unstable weighted approximation (Appendix B), on Nazarov's construction
of spectral measures based on a ``distorted Poisson formula'' (Appendix C),
and on how to get stably orthogonal spectral measures (Appendix D).

\subsubsection*{Acknowledgments} Vladimir Marchenko asked us about the existence
of non-classi\-cal orthogonal spectral functions and explained the relevancy
of the Bernstein weighted approximation problem. Fedor Nazarov generously
helped us with an example of a measure $\mu$ described in the first part of
Theorem~\ref{thm-sharpness} and with an explicit construction of spectral
measures given in Appendix~\ref{appendix-5}. Peter Yuditskii patiently
explained us the main idea of his work~\cite{Yuditskii}. We are grateful to
all of them and to Alexander Ulanovskii for very illuminating discussions.

\section{Weighted exponential approximation on $\R$}\label{se2}

Here, we bring several facts from the well-developed theory of weigh\-ted
exponential approximation on the real axis.

\begin{definition} A lower semi-continuous function
$W\colon\mathbb R\to(0,\infty]$ is called a weight if
$W$ is not equal to $\infty$ identically, and for some $s\in\R$ we have
\begin{equation}
\lim_{|x|\to\infty}(1+|x|)^sW(x)=\infty.
\label{m1}
\end{equation}
\end{definition}

\begin{definition}
$C_0(W)$ is the semi-normed space of the functions $f$ continuous on the
real line and such that
\begin{gather*}
\lim_{|x|\to\infty}\frac{f(x)}{W(x)}=0,\\
\|f\|_W=\sup_{\mathbb R}\frac{|f|}{W}.
\end{gather*}
\end{definition}

\begin{theorem}[M.~Riesz--Mergelyan]\label{la}
Let $\mu$ be a non-negative measure on $\R$ having at most polynomial
growth~\eqref{eq.growth1}, and let $W$ be a weight. Suppose that $a>0$ and
that $X=C_0(W)$ or $X=L^2(\mu)$. The set $\E(a)$ is dense in $X$ if and only
if there exists $f\in\E(a)$ with $f({\rm i})\ne 0$ such
that $f(x)(x-{\rm i})^{-1}\!\in\!\clos_X\E(a)$.
\end{theorem}

This is a counterpart of a classical theorem of M.~Riesz and Mergelyan
(see, for instance, \cite{Mergelyan}). Both of them considered the polynomial approximation,
M.~Riesz in $L^2(\mu)$, while Mergelyan in $C_0(W)$. Rather general versions
of their results can be found in Pitt~\cite[Theorem~3.1]{Pitt} and in
Levin~\cite[Section~I.1]{Levin}.

\medskip\par\noindent{\em Proof of Theorem~\ref{la}:}
The necessity part is evident. To verify the sufficiency part, we fix a
function $f\in\E(a)$ such that $f({\rm i})\ne 0$ and
$f(x)(x-{\rm i})^{-1}\in\clos_X\E(a)$. Without loss of generality, we assume
that $f$ does not vanish on $\R$. Otherwise, we take any function
$g\in\E(a)$ that vanishes at $z=\rm i$ and does not vanish on $\R$, and
replace $f$ by the function $f-cg$ with an appropriate constant $c$
(such a constant exists since the function $f/g$ is analytic in a neighborhood of the real line).
Next, we use that
$$
h\in\E(a) \implies \frac{f({\rm i})h-h({\rm i})f}{\cdot-{\rm i}}\in\E(a).
$$

\begin{claim}\label{claim-RM1} Let $f^*(z)= \overline{f(\overline{z})}$.
Then $f^*(x)(x-{\rm i})^{-1}\in\clos_X\E(a)$.
\end{claim}

\par\noindent{\em Proof:}
We have
\[
\frac{f^*(x)}{x-{\rm i}}
= \frac{f^*({\rm i})}{f({\rm i})}\cdot \frac{f(x)}{x-{\rm i}} +
\frac{f^*(x)f({\rm i})-f(x)f^*({\rm i})}{f({\rm i})(x-{\rm i})}\,.
\]
The first term on the right-hand side belongs to $\clos_X\E(a)$, while
the second term lies in $\E (a)$. We are done. \hfill $\Box$

\begin{claim}\label{claim-RM2}
For each $n\ge 1$ we have
\[
f(x)(x\pm {\rm i})^{-n}\in\clos_X\E(a)\,.
\]
\end{claim}

\par\noindent{\em Proof:}
Let $h_k\in\E(a)$, $h_k\stackrel{X}{\to} f(x)(x-{\rm i})^{-1}$. Then
\[
\big\| \frac{h_k(x)}{x - {\rm i}} - f(x)(x-{\rm i})^{-2} \big\|_X
\le \big\| h_k (x) - f(x)(x - {\rm i})^{-1} \big\|_X \to 0\,,
\]
and the sequence of functions
\[
\frac{h_k(x)}{x-{\rm i}} =
\frac{h_k(x)f({\rm i})
-h_k({\rm i})f(x)}{f({\rm i})(x-{\rm i})}+\frac{h_k({\rm i})f(x)}{f({\rm i})(x-{\rm i})}
\]
is contained in $\clos_X\E(a)$. We obtain that
$f(x)(x-{\rm i})^{-2}\in\clos_X\E(a)$. In the same way,
$f(x)(x-{\rm i})^{-n}\in\clos_X\E(a)$ for each $n\ge 1$.
Furthermore, using Claim~\ref{claim-RM1}, in the same way we obtain that
$f^*(x)(x-{\rm i})^{-n}\in\clos_X\E(a)$ for each $n\ge 1$.
Noting that
\[
f(x)(x+{\rm i})^{-n} = \overline{f^*(x)(x-{\rm i})^{-n}}, \qquad x\in\R\,,
n\ge 1\,,
\]
and that $\clos_X\E(a)$ is closed with respect to the conjugation,
we obtain that $f(x)(x+{\rm i})^{-n}\in\clos_X\E(a)$ for each $n\ge 1$. We are done. \hfill $\Box$

\medskip Let $V_s(x)=(1+|x|)^{-s}$, $C_s(\R)=C_0(V_s)$. For large $s$,
convergence in $C_s(\R)$ implies convergence in $X$.

\begin{claim}\label{claim-RM3}
The linear span of the functions
$\bigl\{f(x)(x\pm{\rm i})^{-n}\bigr\}_{n\ge 1}$ is dense in $C_s(\R)$.
\end{claim}

\par\noindent{\em Proof:}
Otherwise, there is a non-zero finite complex-valued measure $\nu$ on $\R$
such
\[
\int_\R\frac{f(x)(1+|x|)^{s}\,\d\nu(x)}{(x\pm{\rm i})^n}=0\,,\qquad n\in\mathbb N\,,
\]
whence
\[
\int_\R \frac{f(x)(1+|x|)^{s}\, \d\nu (x)}{x - \zeta} = 0\,, \qquad \zeta\in\C\setminus\R\,,
\]
which, in its turn, yields that the measure $f(x)(1+|x|)^{s}\,\d\nu(x)$ vanishes.
Since $f$ does not vanish on
$\R$, $\nu$ is the zero measure, which contradicts our assumption.
\mbox{}\hfill $\Box$

\medskip Now, we easily complete the proof of sufficiency in
Theorem~\ref{la}. By Claim~\ref{claim-RM3}, each continuous function on $\R$
with compact support can be approximated in $C_s(\R)$, and hence in $X$,
by finite linear
combinations of the functions $f(x)(x\pm{\rm i})^{-n}$. Then, by
Claim~\ref{claim-RM2}, continuous functions with compact support belong to
$\clos_X\E(a)$. It remains to recall that continuous functions with compact
support are dense in $X$, completing the proof of Theorem~\ref{la}.
\mbox{} \hfill $\Box$

\medskip Next, we introduce a $C_0(W)$-counterpart of stable density,
cf. Definition~\ref{def-stable_density}.
Given a weight $W$, we set $W_t(x)=W(x)(1+|x|)^{-t}$.
\begin{definition}[stable density in $C_0(W)$]\label{def-stable-C_0(W)}
We say that $\E(a)$ is stably dense in $C_0(W)$ if for each $t\ge 0$,
$\E(a)$ is dense in $C_0(W_t)$.
\end{definition}

The following theorem is a version of a recent result of
Bakan~\cite{Bak} who dealt with weighted polynomial approximation.

\begin{theorem}[Bakan]\label{lb}
Let $\mu$ be a non-negative measure on $\R$ satisfying the growth
condition~\eqref{eq.growth1}, and let $a>0$. The set $\E(a)$ is (stably)
dense in $L^2(\mu)$ if and only if there exists a weight $W\in L^2(\mu)$
satisfying the growth condition~\eqref{m1} such that $\E(a)$ is (stably)
dense in $C_0(W)$.
\end{theorem}

\noindent{\em Proof:} We consider only the stable density case.
The same argument works (with some simplifications) in the other case.

To verify the sufficiency part note that
\begin{multline*}
\|f\|^2_{L^2(\mu_p)}=\int_{\mathbb R}|f(x)|^2 (1+|x|)^p \d\mu(x) \\
=\int_{\mathbb R}\Bigl|\frac{f(x)}{W_{p/2}(x)}\Bigr|^2
W^2(x)\,\d\mu(x)\le \| W\|^2_{L^2(\mu)}\,  \|f\|^2_{W_{p/2}}.
\end{multline*}
Therefore, if there exists a function $f\in\E(a)$ with $f({\rm i})\ne 0$
such that
$$
f(x)(x-{\rm i})^{-1}\in\clos_{C_0(W_{p/2})}\E(a)\,,
$$
then
$$
f(x)(x - {\rm i})^{-1} \in\clos_{L^2(\mu_p)}\E(a)\,.
$$
It remains to apply Theorem~\ref{la}.

To verify the necessity part, we choose $n$ so big that
the function $x\mapsto(1+|x|^n)^{-1}$ belongs to $L^2(\mu)$, and suppose that
$\E(a)$ is dense in $L^2(\mu_p)$
for every $p<\infty$. Take any function $f\in\E(a)$ with $f({\rm i})\ne 0$,
and choose functions $h_k \in\E(a)$ such that
$$
\int_{\mathbb R}|h_k(x)-f(x)(x-{\rm i})^{-1}|^2(1+|x|)^{2k}\,\d\mu(x)<8^{-k}
\,.
$$

We set
\[
W(x)=\Bigl[(1+|x|^n)^{-1}+\sum_{k\ge 1}4^k\bigl|h_k(x)-\frac{f(x)}{x-{\rm i}}\bigr|^2
(1+|x|)^{2k}\Bigr]^{1/2}.
\]
Then $W$ is a lower semi-continuous function (since $W^2$ is the sum
of a series with continuous non-negative terms) satisfying the growth
condition~\eqref{m1}, and $W\in L^2(\mu)$. Let $s\in\R$. For $k\ge s$
we have
\begin{multline*}
\|h_k(x)-f(x)(x-{\rm i})^{-1}\|_{W_s}=\sup_{\mathbb R}\frac{|h_k(x)-
f(x)(x-{\rm i})^{-1}|\,(1+|x|)^s }{W(x)}\\
\le\sup_{\mathbb R}2^{-k}\frac{|h_k(x)-f(x)(x-{\rm i})^{-1}|\,(1+|x|)^s}
{|h_k(x)-f(x)(x-{\rm i})^{-1}|\,(1+|x|)^k}\le 2^{-k}.
\end{multline*}
Thus, by Theorem~\ref{la}, $\E(a)$ is dense in $C_0(W_s)$. \hfill $\square$

\medskip The last result in this section is a version of de~Branges'
classical theorem on weighted polynomial approximation \cite{deBr1},
\cite[Theorem~66]{deBr}, \cite[Section~VI.F]{Koosis1}.
\begin{definition}[entire functions of Krein's class]\label{def-Krein's_class1}\mbox{}
Given $a>0$, we denote by $\mathcal K(a)$ the class of entire
functions $f$ of exponential type $a$ with simple real zeros $\Lambda(f)$
such that $f(\R)\subset\R$, and $1/f$ is represented as an absolutely convergent series:
\begin{equation}
\frac1{f(z)}=R(z)+\sum_{\lambda\in\Lambda(f)}\frac 1{f'(\lambda)}
\Bigl(\frac1{z-\lambda}+\frac1{\lambda}+\frac z{\lambda^2}+\ldots
+\frac {z^{N}}{\lambda^{N+1}}\Bigr)\,,
\label{eq:dop1}
\end{equation}
with some $N\ge 0$ and with a polynomial $R$.

Given a weight $W$, we
denote by $\mathcal K(a,W)$ the class of all functions $f\in \mathcal K(a)$ such that
\begin{equation}\label{eq:dop2}
\sum_{\lambda\in\Lambda(f)}\frac{W(\lambda)}{|f'(\lambda)|}<\infty.
\end{equation}
\end{definition}
By Krein's theorem, entire functions of the class $\mathcal K(a)$ belong to the Cartwright
class~\cite[Theorem~3 in Lecture~16]{Levin-book}. Using the Phragm\'en-Lindel\"of principle,
it is not difficult to verify that if an
entire function $f$ of Cartwright class maps $\R$ into $\R$,
has simple real zeroes $\Lambda (f)$, and satisfies~\eqref{eq:dop2},
then representation~\eqref{eq:dop1} is valid, cf. Kossis~\cite[Section~VIF.4]{Koosis1}.

\begin{theorem}[de Branges]\label{thm-branges} Let $W$ be a weight function,
and let $a>0$. The linear space $\E(a)$ is not
dense in $C_0(W)$ if and only if $\mathcal K(a,W)\ne\emptyset$.
\end{theorem}

\medskip\par\noindent{\em Proof:}
First, we assume that $\mathcal K(a,W)\ne\emptyset$. Let
$B\in\mathcal K(a,W)$. Without loss of generality, assume that $B(0)\not=0$. Then the measure
$$
\mu_B=\sum_{\lambda\in\Lambda(B)}\frac{\delta_\lambda}{B'(\lambda)}
$$
belongs to the dual space $C_0(W)^*$. The Lagrange interpolation
formula shows that, for each $f\in\E(a)$ we have
\[
z f(z) = \sum_{\la\in\Lambda(B)} \frac{\la f(\la) B(z)}{(z-\la) B'(\la)}\,.
\]
Letting $z=0$, we see that the measure $\mu_B$ annihilates $\E(a)$. Hence,
$\E(a)$ is not dense in $C_0(W)$.

\medskip
The other implication is more deep. Here, we
use a modification of the argument presented
in~\cite{SY} for the polynomial approximation problem.
Suppose that $\E(a)$ is not dense in $C_0(W)$.
For $0<b\le a$ denote $X_b = \clos_{C_0(W)}\E(b)$.

Choose $\phi\in\E(a/2)$, real on $\R$, such
that $\phi(\rm i)\not=0$. By Theorem~\ref{la},
the functions $t\mapsto\phi(t)/(t\pm{\rm i})$ do not belong
to $X_b$, $a/2\le b\le a$. The value $a/2$ is of no importance here; 
we can replace it by any number in the interval $(0,a)$.

Given a weight $W$, we consider the semi-normed space $C(W)$ of functions continuous on the real line with finite norm $\|f\|_{W}=\sup_{\,\mathbb R}\,|f|/W$.
Let $Y_b$ be the space of entire functions of exponential type at most $b$ such that 
$f\in C(W)$. Put
$$
M_b(z)=\sup\{|f(z)|\colon f\in\E(b)\cap C(W),\,\|f\|_W\le 1\}.
$$
Note that the function $b\mapsto M_b(z)$ does not decrease.

\begin{claim}\label{claim:1} Let $a/2\le b\le a$. Then
$X_b\subset Y_b$ and  
$$
|f(z)|\le \|f\|_{C(W)}M_b(z),\qquad f\in X_b,\,z\in\mathbb C.
$$
Furthermore,
\begin{gather}
\int_{\R}\frac{\log^{+}M_b(x)}{1+x^2}\, \d x<\infty,\label{m2}\\
\log M_b(z)\le(b+o(1))|z|,\qquad |z|\to\infty\notag.
\end{gather}
\end{claim}

\noindent Proof: 
Let $L$ be a bounded linear functional on $C_0(W)$ vanishing on
$\E(b)$ such that $L\bigl[\phi(t)/(t-{\rm i})\bigr]\not=0$. The
function $F\colon w\mapsto L \bigl[ \phi/(\cdot-w) \bigr]$ is analytic
in the upper half-plane, and
$$
\int_{\R}\frac{\log^{-}|F(x+{\rm i})|}{1+x^2}\, \d x<\infty.
$$
Then the relation
$$
f(w)L\Bigl[ \frac{\phi}{\cdot-w}\Bigr] = \phi (w) L\Bigl[ \frac{f}{\cdot-w}
\Bigr],\qquad f\in\E(b)\cap C(W),
$$
yields that
$$
\int_{\R}\frac{\log^{+}M_b(x+{\rm i})}{1+x^2}\, \d x<\infty.
$$
Then estimates for subharmonic functions in a half plane imply \eqref{m2}. Furthermore, for every $\varepsilon>0$,
$$
\log M_b(w)\le(b+o(1))|\Im w|,\qquad\epsilon<\arg w<\pi-\epsilon,\,|w|\to
\infty.
$$
A Phragm\'en-Lindel\"of type theorem completes our argument; see
also~\cite[VI~E]{Koosis1}.
\hbox{} \hfill $\Box$

\begin{claim}\label{claim:2} Let $a/2\le b<a$. Then $Y_b$ is a closed subset of 
$C(W)$ and 
\begin{equation}
|f(z)|\le \|f\|_{C(W)}\inf_{b<d\le a}M_d(z),\qquad f\in Y_b,\,z\in\mathbb C.
\label{192}
\end{equation}
\end{claim}

\noindent Proof:  Let $f\in Y_b$, and let 
$0<\varepsilon<(a-b)/3$. The function 
$$
x\mapsto f(x)\sin^2(\varepsilon x)/(\varepsilon x)^2
$$
is an entire function of exponential type less than $b+3\varepsilon$ belonging to
$C_0(W)$. It is a classical fact that such a function belongs to 
$X_{b+3\varepsilon}$, cf. \cite[Section~VI~H~1]{Koosis1}, \cite[Section 2]{LY}.
This gives \eqref{192}.

Furthermore, if $f_n\in Y_b$, $n\ge 1$, and the sequence $(f_n)_{n\ge 1}$ converges in $C(W)$, then 
by \eqref{192}, $(f_n)_{n\ge 1}$ converges uniformly on the compact subsets of the complex plane
to an entire function of exponential type at most $b$ belonging to $C(W)$.
\hbox{} \hfill $\Box$
\smallskip

Next, consider the linear space 
$$
Z_b=\{(1+x^2)f:f\in Y_b\}.
$$
Since $Y_b$ is a closed subset of $C(W)$, we obtain that 
$Z_b$ is closed in $C(W_{-2})$ (as above, $W_t(x)=W(x)(1+|x|)^{-t}$).
Furthermore, since
$$
\frac{\phi(t)(t\pm{\rm i})}{t^2+1}
=\frac{\phi(t)}{t\mp{\rm i}},
$$
we can find a function $\psi\in\E(a/2)$ ($\psi(t)=\phi(t)$
or $\psi(t)=t\phi(t)$), real on $\R$, such
that $\psi(\rm i)\not=0$, and $\widetilde{\psi}\not\in X_a$, 
$\widetilde{\psi}(t)=\psi(t)/(t^2+1)$.

For $a/2\le b<a$ consider a Chebyshev type extremal problem of the best
approximation to $\psi$ by elements of $Z_b$ in the
$\|\cdot\|_{W_{-2}}$-norm.
By a normal family argument, there exists a function $f_b\in Z_b$
such that 
$$
\|f_b-\psi\|_{W_{-2}}=\dist_{C(W_{-2})}(\psi,Z_b)=L_b\ge 0.
$$
We can assume that $f_b$ is real on the real line.
Note that the function $b\mapsto L_b$ does not increase.

For $0<\varepsilon<(a-b)/2$, the functions
$$
f_{b,\varepsilon}(x)=\frac{f_b(x)}{x^2+1}\cdot \frac{\sin^2(\varepsilon x)}{(\varepsilon x)^2}
$$
belong to $X_a$, and hence,
$$
L_b=\|f_b-\psi\|_{W_{-2}}=\lim_{\varepsilon\to 0+}\|(x^2+1)f_{b,\varepsilon}-\psi\|_{W_{-2}}
\ge c\,\dist_{C(W)}(\widetilde{\psi},X_a),\quad a/2\le b<a.
$$
Therefore,
$$
L=\lim_{b\to a-0}L_b>0.
$$

Put $g_b=\psi-f_b$.

\begin{claim}\label{claim:3} The function $g_b$ is of exponential type $b$, has only
simple real zeros, and between each pair of consecutive zeros there is a
point $\lambda$ such that $|g_b(\lambda)|=L_b\,W_{-2}(\lambda)$.
\end{claim}

\noindent Proof: uses so called ``Markov's corrections''.
If one of the assertions of the lemma does not hold, then we construct
an entire function
\[
f_*(z)  = f_b(z) + (z^2+1)g_b(z)R(z)
\]
such that $g_bR\in Y_b$, and
$$
\sup_\R \bigl|1 - (x^2+1)\, R(x) \bigr|\, \frac{|g_b(x)|}{W_{-2}(x)} <
\sup_\R \frac{|g_b|}{W_{-2}}\,. 
$$
Then $f_*\in Z_b$. Since $ \psi - f_*= (\psi - f_b) (1-(x^2+1)\, R)$, 
we will conclude that $\| \psi - f_*\|_{W_{-2}} < \| \psi - f_b\|_{W_{-2}}$, which 
contradicts to the extremality property of the function $f_b$.
\medskip

First, we show that $g_b$ is of exponential type $b$. Otherwise,
considering the functions $R$,
$$
R(x)=\frac12\Bigl[\frac{\sin^2(\epsilon x)}{x^2}
+\frac{\sin^2(\epsilon (x+1))}{(x+1)^2}\Bigr], \qquad\epsilon>0,
$$
for small positive $\epsilon$,
we obtain that $f_b$ is not extremal. Indeed,
$$
\sup_\R\, \Bigl|1-\frac{x^2+1}2
\Bigl[\frac{\sin^2(\epsilon x)}{x^2}
+\frac{\sin^2(\epsilon (x+1))}{(x+1)^2}\Bigr]\Bigr|<1,
$$
and the functions
$$
x\mapsto g_b(x)
\Bigl[\frac{\sin^2(\epsilon x)}{x^2}
+\frac{\sin^2(\epsilon (x+1))}{(x+1)^2}\Bigr]
$$
belong to $Y_b$, for small $\epsilon>0$.

\medskip
Next, $g_b$ has only simple real zeros. Otherwise, if $\alpha,\bar\alpha$ is a
pair of conjugate zeros of $g_b$, then we set 
$$
R(x)=\frac{\epsilon}{(x-\alpha)(x-\bar\alpha)}, \qquad \epsilon>0,
$$
and use that $g_bR\in Y_b$.

\medskip
Finally, if $\alpha<\beta$ are two consecutive zeros of $g_b$ such that $|g_b|<L_b\,W_{-2}$
on $(\alpha,\beta)$, then by the lower semicontinuity of $W_{-2}$, we have
$\sup_{(\alpha,\beta)}\,|g_b|/ W_{-2}<L_b$. It remains to consider the functions $R$,
$$
R(x)=\frac{\epsilon}{(x-\alpha)(x-\beta)} 
$$
for small positive $\epsilon>0$, to show that $g_b$ is not extremal. 
\hfill $\square$
\medskip

Consider the sequence of entire function $(g_{a-\frac1n})_{n\ge 1}$.
By Claim~\ref{claim:2}, it is a normal family, and we can find a limit function $g$.
Then $g$ is of exponential type at most $a$,
$$
\|g\|_{C(W_{-2})}\le L=\lim_{b\to a-0}L_b,
$$
and $g({\rm i})=\lim_{b\to a-0}g_b({\rm i})=\psi({\rm i})\not=0$.

Set $f=\psi-g$. If $g$ (and, hence, $f$) is of exponential type $b<a$, then, arguing as in the first part of the proof of 
Claim~\ref{claim:3}, and comparing $f$ to $f_{\max(b,a/2)}$, we get a contradiction  
to the extremality property of the function $f_{\max(b,a/2)}$.
Thus, $g$ is of exponential type $a$; it is real on the real line and its zeros are real. 
Using the lower semicontinuity of $W_{-2}$, we obtain that the zeros of $g$ are simple, and
between each pair of consecutive zeros of $g$ there is a
point $\lambda$ such that $|g(\lambda)|=L\,W_{-2}(\lambda)$. In particular, this implies that $\|g\|_{C(W_{-2})}=L$.

\medskip
Now we complete the proof of the de Branges theorem. Let
$\Lambda=\Lambda(g)=\{x_j\}$, and let
$\Lambda^*=\{x'_j\}\subset\{x:|g(x)|=L\,W_{-2}(x)\}$, with interlacing $x_j$ and
$x'_j$: $x_j<x'_j<x_{j+1}$. Set
$$
B(z)=g(z)\prod_j\frac{z-x'_j}{z-x_j}.
$$
Then $\Im(g/B)$ does not change the sign in the upper half-plane, and hence,
$$
\sum_{\lambda\in\Lambda^*}\frac{|g(\lambda)|}{|B'(\lambda)|(1+\lambda^2)}<\infty.
$$
By the definition of $\Lambda^*$,
$$
\sum_{\lambda\in\Lambda^*}\frac{W(\lambda)}{|B'(\lambda)|}<\infty.
$$
Finally, by Claim~\ref{claim:2}, 
$$
|g_b(x)|\le CM_a(x),\qquad a/2\le b<a,
$$
and by Claim~\ref{claim:1} we obtain that $g$, and hence, $B$ 
belong to the Cartwright class. We conclude that
$B\in\mathcal K(a,W)$.
\mbox{}\hfill $\square$

\section{Proof of Theorem~\ref{thm_main1}}\label{se3}

First, let us recall the definition \eqref{m0} of the intervals $I_x$,
$kI_x$. An elementary calculation shows that for $k_1=k_1(\delta)$,
we have the following property:
$$
y\in I_x\implies 2I_x\subset k_1I_y.
$$

\begin{lemma} \label{le}
Suppose that $\mu \preccurlyeq \widetilde\mu$. Given $p>0$ and a weight
function $\widetilde W\in L^2(\widetilde\mu)$, we define a function $W$ by
\begin{equation}
W(x)=\min\bigl[\inf_{k_1I_x}\widetilde W_p, e^{\delta|x|/3}\bigr].
\label{m4}
\end{equation}
Then $W$ is a weight function, and $W\in L^2(\mu)$, provided that $p$ is big
enough.

Here, as above, $\widetilde W_p(x) = (1+|x|)^{-p}\widetilde W(x)$.
\end{lemma}

\noindent{\em Proof:} It is immediately seen that $W$ is a weight function.
Next we choose $x_j$, $j\in\mathbb Z$, on $\R$ in such a way that the
intervals $I_{x_j}$ cover $\R$ with intersections only at endpoints. Then
the intervals $2I_{x_j}$ cover $\R$ with multiplicity bounded by
$C=C(\delta)$. Using that
\begin{gather*}
\sup_{y\in I_x}W(y)\le \sup_{y\in I_x}\inf_{t\in k_1I_y}\widetilde W_p(t)
\le \sup_{y\in I_x}\inf_{t\in 2I_x}\widetilde W_p(t)
=\inf_{t\in 2I_x}\widetilde W_p(t),\\
\mu(I_x)\le C(1+|x|)^n\left(\widetilde\mu(2I_x)+e^{-2\delta|x|} \right),
\end{gather*}
we obtain that for $p\ge n$,
\begin{multline*}
\int_{I_x}W^2(x)\,d\mu(x)\le \sup_{I_x}W^2\cdot\mu(I_x)\\
\le C(1+|x|)^{2n} \widetilde\mu(2I_x)\cdot\inf_{2I_x}\widetilde W_p^2+
C (1+|x|)^{2n} e^{-2\delta|x|}\cdot e^{2\delta|x|/3}\\
\le C (1+|x|)^{2n-2p} \int_{2I_x}\widetilde W^2(y)\, \d\widetilde \mu(y)+
C  \int_{2I_x} (1+|y|)^{2n}e^{-\delta|y|/3}\,\d y.
\end{multline*}
Summing up these inequalities for $x=x_j$, $j\in\mathbb Z$, we complete the proof. \hfill $\square$

\begin{lemma} \label{lf} Let $W$ be a weight, and let $B\in\mathcal K(a,W)$.
Then for some $c>0$ and $C<\infty$ we have
\begin{gather}
|B(x)|+|B'(x)|+|B''(x)|<e^{\delta|x|/5},\qquad |x|>C,\label{m5}\\
|\lambda-\lambda'|>c e^{-\delta|\la|/4},\qquad \lambda,\lambda'\in
\Lambda(B),\,\lambda\ne\lambda',\label{m6}\\
\sum_{\lambda\in\Lambda(B),\,\lambda\ne 0}\frac{1}{|\lambda|^2}<\infty.
\label{m7}
\end{gather}
\end{lemma}

\noindent{\em Proof:} Since $B$ is of exponential type, we have \eqref{m7}.
Denote by $H_B$ the Phrag\-m\'en-Lindel\"of indicator function of the entire
function $B$,
$$
H_B(\theta)=\limsup_{r\to\infty} \frac{\log |B(re^{i\theta})|}{r}.
$$
Since the function $B$ belongs to the Cartwright class, we have
$H_B(0)=H_B(\pi)=0$ (see \cite[Lecture~16]{Levin-book}). As a consequence of
Cauchy's formula for the derivative, the indicator of the derivative of an
entire function cannot exceed the indicator of the function itself. Hence,
$$
H_{B'}(0)\le 0,\quad H_{B'}(\pi)\le 0,\quad H_{B''}(0)\le 0,\quad
H_{B''}(\pi)\le 0,
$$
which implies \eqref{m5}.

Using that
$$
\sum_{\lambda\in\Lambda(B)}\frac{W(\lambda)}{|B'(\lambda)|}<\infty
$$
and that
$$
\lim_{|x|\to\infty}(1+|x|)^s W(x)=\infty
$$
for some $s\in\R$, we obtain
$$
|B'(\lambda)|\ge \frac {c}{(1+|\lambda|)^s},\qquad \lambda\in\Lambda(B).
$$
Since the signs of $B'$ at consecutive zeros of $B$ are opposite, this
inequality together with estimate \eqref{m5} on $B''$ gives \eqref{m6}. \hfill
$\square$

\begin{lemma} \label{lg} Let $\widetilde W$ be a weight. If
$\mathcal K(a,W)\not=\emptyset$, where $W$ is defined by \eqref{m4},
then $\mathcal K(a,{\widetilde W}_{p+\ell})\not=\emptyset$ provided that
$\ell$ is big enough.
\end{lemma}

\noindent{\em Proof:} Let $B\in\mathcal K(a,W)$. Since
$$
\sum_{\lambda\in\Lambda(B)}\frac{W(\lambda)}{|B'(\lambda)|}<\infty,
$$
for some $c>0$ we have
$$
|B'(\lambda)|\ge c\cdot\min\bigl[\inf_{k_1I_\lambda}\widetilde W_p,
e^{\delta|\lambda|/3}\bigr],\qquad \lambda\in \Lambda(B).
$$
By \eqref{m5},
$$
|B'(\lambda)|=o(e^{\delta|\lambda|/3}),\qquad |\lambda|\to\infty,
$$
and, hence, for some $C<\infty$,
$$
|B'(\lambda)|\ge c\cdot\inf_{k_1I_\lambda}\widetilde W_p,
\qquad \lambda\in \Lambda(B),\,|\lambda|>C.
$$
Let $D_\lambda$ be the disc centered at $\lambda\in\R$ of radius
$e^{-\delta|\lambda|/3}$. For some $M<\infty$, by \eqref{m6} we have the
following implication:
$$
\lambda,\lambda'\in\Lambda(B),\quad |\lambda|>M,\,\,|\lambda'|>M, \quad
\lambda\not=\lambda' \implies D_\lambda\cap D_{\lambda'}=\emptyset.
$$
and $k_1I_{\lambda}\subset D_\lambda$ for $\lambda\in\Lambda(B)$,
$|\lambda|>M$.

Now, for some $c>0$ and for every $\lambda\in\Lambda(B)$ with $|\lambda|>M$
we find $\zeta_\lambda\in k_1I_\lambda$ such that
$$
|B'(\lambda)|\ge c\widetilde W_p(\zeta_\lambda).
$$

Without loss of generality assume that $B(0)=1$. Since $B$ is of Cartwright
class, we have
$$
B(z)=\lim_{R\to\infty}\prod_{|\lambda| \le R,\,\lambda\in\Lambda(B)}
\Bigl(1-\frac{z}{\lambda} \Bigr).
$$

For $z\in \partial D_\lambda$, $\lambda\in\Lambda(B)$, $|\lambda|>M$ we have
\begin{equation}
\Bigl|\Bigl(1-\frac z{\zeta_\lambda}\Bigr)
\Bigl(1-\frac z{\lambda}\Bigr)^{-1}-1\Bigr|=
\Bigl|\frac z{\zeta_\lambda}\cdot
\frac {\lambda-\zeta_\lambda}{z-\lambda}\Bigr|
\le ce^{-2\delta|\lambda|/3},
\label{xx1}
\end{equation}
and by the maximum principle, we have the same estimate for all
$z\in\mathbb C\setminus D_\lambda$.

We define an entire function $B_1$ by
$$
B_1(z)=\lim_{R\to\infty}
\prod_{M < |\lambda| \le R,\,\lambda\in\Lambda(B)}
\Bigl(1-\frac{z}{\zeta_\lambda} \Bigr);
$$
the limit on the right hand side exists because of \eqref{m7} and
\eqref{xx1}. In a similar way, applying the maximum principle in every
$D_\lambda$, we conclude that $B_1$ is of exponential type. Furthermore,
$$
|B_1(z)|\ge c|B(z)|(1+|z|)^{-N},\qquad z\in \partial D_\lambda,\,
\lambda\in\Lambda(B)\setminus [-M, M],
$$
where $N=\card(\Lambda(B)\cap[-M,M])$, and, hence,
$$
|B'_1(\zeta_\lambda)|\ge c (1+|\zeta_\lambda|)^{-N} |B'(\lambda)|,
\qquad \lambda\in\Lambda(B)\setminus [-M, M],
$$
which implies that $B_1\in\mathcal K(a,{\widetilde W}_{p+\ell})$ with
$\ell=N+2$. \hfill $\square$

\medskip Now we are ready to pass to

\medskip\par\noindent{\em Proof of Theorem~\ref{thm_main1}}:
Let $\E(a)$ be stably dense in $L^2(\widetilde\mu)$. By Bakan's
Theorem~\ref{lb}, there exists a weight $\widetilde W\in L^2(\widetilde\mu)$
such that $\E(a)$ is stably dense in $C_0(\widetilde W)$. Then, by de
Branges' Theorem~\ref{thm-branges}, for each $t<\infty$, we have
$\mathcal K(a,\widetilde W_t)=\emptyset$. We take $p$ big enough. Then the
function $W$ defined in \eqref{m4} is a weight function and belongs to
$L^2(\mu)$ by Lemma~\ref{le}. Hence, $\mathcal K(a,W)=\emptyset$ (otherwise,
by Lemma~\ref{lg}, $\mathcal K(a,\widetilde W_{p+\ell})\ne\emptyset$ for large $p$, which
is impossible). Now, using again de Branges' Theorem, we obtain that $\E(a)$
is dense in $C_0(W)$. Applying again Bakan's Theorem, we see that $\E(a)$ is
dense in $L^2(\mu)$, proving the theorem. \mbox{}\hfill $\square$

\subsection{Remark on the polynomial approximation in $L^2(\mu)$} \label{add1}

Theorem~\ref{thm_main1} has a polynomial counterpart, which can be proved
using the same lines of reasoning.
{\em Let $\widetilde\mu$ be a non-negative measure
with finite moments,
$$
\int_{\mathbb R}|x|^n\, {\rm d} x<\infty,\qquad n\ge 0,
$$
and let $\mu\preccurlyeq\widetilde\mu$. If the set of the polynomials
$\mathcal P$ is stably
dense\footnote{That is, $\mathcal P$ is dense in
$L^2(\widetilde{\mu}_t)$ for each $t<\infty$.
Equivalently, one can say that the measure $\widetilde{\mu}$ has
infinite index of determinacy
for the Hamburger moment problem on the real axis, see Berg and
Duran~\cite{BergDuran}}
in $L^2(\widetilde{\mu})$, then $\mathcal P$ is dense in $L^2(\mu)$}.
In other words,
{\em if the measure $\widetilde{\mu}$ has infinite index of
determinacy for the Hamburger moment
problem, and $\mu\preccurlyeq \widetilde{\mu}$, then the measure $\mu$
is determinate}.

This extends a result of Yuditskii~\cite{Yuditskii}. Answering a
question posed by Berg in the 1990-s,
Yuditskii obtained the same conclusion with a
much stronger assumption. He proved that if a measure $\widetilde{\mu}$
has infinite index
of determinacy and $\mu = \widetilde{\mu} + \nu$ where $\nu$ is a non-negative
measure on $\mathbb R$ with a finite exponential moment: for some $\delta>0$,
\[
\int_{\mathbb R} e^{\delta|\lambda|}\, {\rm d}\nu (\lambda) > 0\,,
\]
then the measure $\mu$ is determinate.
Peter Yuditskii showed us that (after an appropriate
adjustment) this also
holds for continuous analogues of the Hamburger moment problem
(see~\cite[Chapter~V, \S~3]{Akh}
and items~10--12 in Addenda and Problems to this Chapter) that
correspond to the exponential
approximation problem and to the inverse spectral theory of the one-dimensional
Schr\"odinger equation.

It would be interesting to check whether the elegant
operator-theoretical approach developed by
Yuditskii in~\cite{Yuditskii} (see also Section~\ref{add2} below)
allows one to treat general perturbations
of the kind $\mu\preccurlyeq \widetilde{\mu}$.

\section{Sharpness. Proof of Theorem~\ref{thm-sharpness}}\label{se4}

\subsection{Proof of Theorem~\ref{thm-sharpness}~(i)}

We construct here a convex function $f$ on $[0,\infty)$ such that
\begin{gather}
\int^\infty f(t)e^{-t}\,\d t=\infty,\label{xx5}\\
\int^\infty \min\bigl(f(t),\varepsilon(e^{t})e^{t}\bigr)e^{-t}\,\d t<\infty,\label{xx6}
\end{gather}
and set $\varphi(x)=e^{-f(\log^+|x|)}$, $x\in\mathbb R$.

Then
$$
\int_{\mathbb R}\frac{\log(1/\varphi(x))}{x^2+1}\,\d x=\infty,
$$
and by the result mentioned in Appendix~\ref{xx2}, we have
$T(\varphi(x)\,\d x)=0$. (In fact, a classical theorem of Izumi--Kawata (see
\cite[VID, p.170]{Koosis1}) shows that already the polynomials are dense in
$L^2(\varphi(x)\,dx)$.)

On the other hand, if $\psi(x)=\varphi(x)+e^{-\varepsilon(|x|)|x|}$, then
$$
\int_{\mathbb R}\frac{\log(1/\psi(x))}{x^2+1}\,\d x<\infty,
$$
and by the Krein theorem mentioned in Appendix~\ref{xx3},
$T(\psi(x)\,dx)=\infty$.

The function $f$ will be built as the sum of the functions $f_n=\max(l_n,0)$ for
some linear functions $l_n$ chosen in an inductive process. On step $n\ge 1$
we fix a sufficiently large $a$ such that
$$
\gamma=\sup_{t\ge a}\varepsilon(e^t)<4^{-n},
$$
and $b>a+1$ such that
\begin{gather}
(b-a)\gamma<2^{-n},\label{xx4}\\
\frac{e^{b-a}-1}{b-a}\gamma\ge 10.\label{xx4a}
\end{gather}
The linear function $l_n$ is determined by the conditions
$$
l_n(a)=\gamma e^a,\qquad l_n(b)=\gamma e^b.
$$
To prove \eqref{xx5}--\eqref{xx6} we need only to verify that
\begin{gather}
\int_0^af_n(x)e^{-x}\,\d x+\int_b^{\infty}f_n(x)e^{-x}\,\d x\le C\gamma,
\notag
\\
\int^b_af_n(x)e^{-x}\,\d x\ge 1,\notag
\\
\int^b_a\gamma e^x\cdot e^{-x}\,\d x<2^{-n}.\label{xx10}
\end{gather}
By \eqref{xx4a}, we have
$$
\frac{e^b-e^a}{b-a}\ge \frac{10e^a}{\gamma}\ge e^a.
$$
Therefore, $f_n=0$ on $[0,a-1]$, and
$$
\int_0^af_n(x)e^{-x}\,\d x\le
\gamma\int_{a-1}^a e^{a-x}\,\d x= (e-1)\gamma.
$$
Furthermore,
\begin{gather*}
\int_b^\infty f_n(x)e^{-x}\,\d x= \gamma\int_0^\infty
\Bigl(e^b+\frac{e^b-e^a}{b-a}s\Bigr)e^{-b-s}\,\d s\\=
\gamma\int_0^\infty \Bigl(e^{-s}+\frac{1-e^{a-b}}{b-a}se^{-s}\Bigr)\,\d s\le 2\gamma,
\end{gather*}
and by \eqref{xx4a},
\begin{gather*}
\int_a^b f_n(x)e^{-x}\,\d x= \gamma\int_0^{b-a}
\Bigl(e^a+\frac{e^b-e^a}{b-a}s\Bigr)e^{-a-s}\,\d s\\=
\gamma\int_0^{b-a} \Bigl(e^{-s}+\frac{e^{b-a}-1}{b-a}se^{-s}\Bigr)\,\d s\ge
\gamma\frac{e^{b-a}-1}{b-a}\int_0^1 se^{-s}\,\d s
\ge 1.
\end{gather*}
Finally, \eqref{xx10} follows from \eqref{xx4}. \hfill $\square$
\bigskip

\subsection{Proof of Theorem~\ref{thm-sharpness}~(ii)}

Here, the construction is more involved.

Without loss of generality, we can assume that $\epsilon$ does not increase.

\begin{lemma} \label{lq1} Let $\epsilon(r)\searrow 0$, $r\to\infty$. There
exists a system of disjoint intervals $I_k=[y_k,2y_k]$, $k\ge 1$, and a
convex function $\varphi$ on $[1,\infty)$ such that
\begin{equation}
\epsilon(e^x)e^x=o(\varphi(x)),\qquad e^x\in \cup_{k\ge 1} I_k,\,x\to+\infty,
\label{q2}
\end{equation}
and
\begin{equation}
\int_0^\infty \varphi(t)e^{-t}\, dt<\infty.
\label{q3}
\end{equation}
\end{lemma}

\noindent{\em Proof}: We construct the function $\varphi$ as the sum of functions
$\varphi_k$,
$$
\varphi_k(x)=\max\bigl(\gamma_k y_k(x+1-\log y_k),0\bigr),
$$
with
$$
\gamma_k=k\sup_{t\ge y_k}\epsilon(t).
$$
Then \eqref{q2} follows immediately. Since
\begin{multline*}
\int_0^\infty \varphi_k(t)e^{-t}\, dt=\gamma_k y_k \int_{(\log y_k)-1}^\infty (t+1-\log y_k)e^{-t}\, dt\\=
e\gamma_k\int_0^\infty te^{-t}\, dt=e\gamma_k,
\end{multline*}
we can find a sequence $\{y_k\}$ such that $I_k$ are disjoint,
$\sum_{k\ge 1}\gamma_k<\infty$, and hence, \eqref{q3} holds.\hfill $\square$
\medskip

By Lemma~\ref{lq1}, we obtain $\varphi$ and $\{I_k\}$, and introduce
an even weight
\begin{equation}
W(x)=\exp \varphi(\max(\log|2x|,1)).
\label{q4}
\end{equation}

\begin{definition}
\label{def-Hamb_class1}
We denote by $\mathcal H(W)$ the Hamburger class of transcendental entire
functions $f$ of zero exponential type with simple real zeros
$\Lambda(f)$ such that $f(\R)\subset \R$, and
$$
\sum_{\la\in\Lambda(f)}\frac{W(\la)}{|f'(\la)|}<\infty\,.
$$
\end{definition}

\begin{lemma} \label{lq5} There exists $F\in\mathcal H(W)$ such that for
some $c>0$ and $E\subset\mathbb R$ of finite length, symmetric with respect
to $0$, we have
\begin{gather}
|F(x)|\ge cW(x/2)^c,\qquad x\in\mathbb R\setminus E,
\label{q9}\\
\dist(\lambda,\mathbb R\setminus E)\ge \frac{c}{1+|\lambda|^2},\qquad \lambda\in\Lambda(F).
\label{q10}
\end{gather}
\end{lemma}

\noindent{\em Proof:} We use here an idea from \cite{Sodin}.
{\bf A.} First, we check that there is an entire function $F$ in
$\mathcal H(W)$ with the zero set $\Lambda_F$ symmetric with respect to
the origin. This will readily follow from a version of de~Branges' theorem
dealing with weighted polynomial approximation. We have
$$
\lim_{x\to\infty}\frac{\log W(x)}{\log x}=\infty,
$$
and the polynomials belong to $C_0(W)$. Furthermore, by \eqref{q3},
$$
\int_{-\infty}^{\infty} \frac{\log W(x)}{1+x^2}\,dx<\infty,
$$
and by the Hall theorem \cite[VI\,D]{Koosis1}, the polynomials are not dense
in $C_0(W)$.

Let $\mu\in(C_0(W))^*$, $\mu\not=0$, vanish on the polynomials. Consider the
functional $\tilde\mu\in(C_0(W))^*$ defined by
$$
\langle \tilde\mu,f\rangle =\langle \mu, x\mapsto f(-x)\rangle,
$$
and put $\mu_{\text{even}}=(\mu+\tilde\mu)/2$,
$\mu_{\text{odd}}=(\mu-\tilde\mu)/2$. Suppose that $\mu_{\text{even}}\not=0$
(the case $\mu_{\text{odd}}\not=0$ is dealt with analogously).
Set $W_0(x)=W(\sqrt{x})$ if $x\ge 0$, $W_0(x)=\infty$ otherwise,
$$
\langle \mu_{\text{right}},f\rangle =\langle \mu_{\text{even}}, x\mapsto f(x^2)\rangle.
$$
Then $\mu_{\text{right}}\in(C_0(W_0))^*$, $\mu_{\text{right}}\not=0$, and
$\mu_{\text{right}}$ vanishes on the polynomials. By the de Branges theorem \cite[VI\,F2, VI\,F1]{Koosis1}, there exists a transcendental entire function $F_0$ of at most minimal type of order $1/2$,
$$
\limsup_{|z|\to\infty}\frac{\log|F_0(z)|}{|z|^{1/2}}=0,
$$
real on the real line, with zeros $\{x_k^2\}_{k\ge 1}$, $x_k\ge 0$, $k=o(x_k)$, $k\to\infty$, such that
$$
\sum_{k\ge 1}\frac{W_0(x_k^2)}{|F_0'(x_k^2)|}<\infty.
$$

Let
\begin{equation}
F(z)=F_0(z^2)=F(0)\prod_{k\ge 1}\Bigl(1-\frac{z^2}{x_k^2}\Bigr)
\label{q6}
\end{equation}
(with an obvious modification if $x_1=0$). Then
$$
\sum_{k\ge 1}\frac{x_kW(x_k)}{|F'(x_k)|}<\infty,
$$
and $F\in \mathcal H(W)$.

\smallskip
{\bf B.} Now, we prove estimates~\eqref{q9} and \eqref{q10}. Suppose that
$0<x_k<x_{k+1}$ are two consecutive zeros of $F$. Suppose that
$\Delta=x_{k+1}-x_k>k^{-2}$ (otherwise, just add the closure of the interval
$J_k=(x_k,x_{k+1})$ to $E$). By the Laguerre theorem, the zeros of $F$ and
$F'$ interlace. Denote by $\lambda$ the zero of $F'$ on $J_k$, and set
$G(z)=F'(z)/(z-\lambda)$. Then $G$ has no zeros in the strip
$J_k+i\mathbb R$. Since $|G(x+iy)|$ increases in $y$ for positive $y$, we
have
\begin{gather*}
|G(x_k+{\rm i}y)| \ge |G(x_k)|\ge c W(x_k)/\Delta,  \\
|G(x_{k+1}+{\rm i}y)|\ge |G(x_{k+1})|\ge c W(x_{k+1})/\Delta
\end{gather*}
for $y\in\R$. Hence, by the three lines theorem applied in the strip
$J_k+i\mathbb R$ to the harmonic function $-\log|G(z)|$, we obtain
$$
\log|G(x)|\ge \frac {x-x_k}{\Delta}\log W(x_{k+1})+\frac {x_{k+1}-x}{\Delta}\log W(x_{k})+\log\frac{c}{\Delta},
\,\, x\in J_k.
$$

If $x_{k+1}<2x_k$, then we obtain
\begin{gather*}
\log|G(x)|\ge c\log W(x/2),\qquad x\in J_k,\\
|F'(x)|\ge c W(x/2)^c, \qquad x\in J_k,\,|x-\lambda|>k^{-2},\\
|F(x)|>c_1W(x/2)^{c_1},\qquad  x_k+k^{-2}<x<x_{k+1}-k^{-2},\, |x-\lambda|>k^{-2}.
\end{gather*}

Otherwise, if $x_{k+1}\ge 2x_k$, then, in the same way, for some
$x_k',x_{k+1}'$ with $x_k<x_k'<x_k+k^{-2}$,
$x_{k+1}-k^{-2}<x_{k+1}'<x_{k+1}$, we obtain:
\begin{gather*}
|F(x_k')|\ge cW(x_k/2)^{c},\\
|F(x_{k+1}')|\ge cW(x_{k+1}/2)^{c}.
\end{gather*}
Since the function $|F|$ is log-concave on $J_k$ (this follows immediately
from the representation \eqref{q6}), and $W$ is log-convex, we conclude that
$$
|F(x)|\ge cW(x/2)^{c},\qquad  x_k'\le x\le x_{k+1}',
$$
proving the lemma. \hfill $\square$

\smallskip
Next we use the following simple lemma on perturbations of the sine
function.

\begin{lemma} \label{lq7} Let $\Sigma=\bigcup_{k\ge 0}\{a_k,b_k,c_k,d_k\}$, where for
each $k\ge 0$,
\begin{gather*}
2k+\frac 65<a_k<b_k<2k+\frac 75, \quad 2k+\frac 85<c_k<d_k<2k+\frac 95, \\
\eta_k=d_k-c_k=b_k-a_k, \quad a_k+b_k+c_k+d_k=8k+6\,.
\end{gather*}
Put
$$
G(z)=z\prod_{\lambda\in\Sigma}\Bigl(1-\frac{z^2}{\lambda^2}\Bigr).
$$
Then $G$ is an entire function of exponential type $2\pi$, and
\begin{equation}
|G'(\lambda)|\ge c\eta_k,\qquad \lambda\in \Sigma\cap(2k+1,2k+2).
\label{q8}
\end{equation}
\end{lemma}

\noindent{\em Proof:} (compare to that of Lemma~\ref{lg}) Let $D_k$ be the disc centered at $2k+\frac 32$
of radius $2/3$,
\begin{gather*}
g_k(z)=\frac{(1-z/a_k)(1-z/b_k)(1-z/c_k)(1-z/d_k)}{(1-z/(2k+1))^2(1-z/(2k+2))^2},\\
b_k=\frac{a_kb_kc_kd_k}{(2k+1)^2(2k+2)^2}.
\end{gather*}
Then
\begin{gather*}
b_k=1+O(1/k^2),\qquad k\to\infty,\\
|b_kg_k(z)-1|\le \frac{M}{1+|z-2k|^2}, \qquad z\in\mathbb C\setminus D_k,
\end{gather*}
with $M$ independent of $k$. Therefore,
$$
G_0(z)=\prod_{k\ge 0}[g_k(z)g_k(-z)]
$$
is bounded outside $\bigcup_{k\ge 0}(D_k\cup\widetilde{D_k})$, where
$\widetilde{D_k}=\{w:-w\in D_k\}$, and
\begin{gather*}
|G_0(z)| \asymp |g_k(z)|,\qquad z\in D_k,\\
|G_0(z)| \asymp |g_k(-z)|,\qquad z\in \widetilde {D_k}.
\end{gather*}
Using the maximum principle in the discs $D_k$ and $\widetilde{D_k}$, we
conclude that
$$
G(z)=G_0(z)\cdot \sin^2 \pi z
$$
is an entire function of exponential type $2\pi$; estimate \eqref{q8}
follows immediately. \mbox{}\hfill $\square$

\smallskip
Now we return to the proof of Theorem~\ref{thm-sharpness}~(ii). Let the function $F$ and the set $E$ be as in Lemma~\ref{lq5}. Put
\begin{align*}
A=&\{k\ge 0:(2k+1,2k+2)\not\subset\cup_{n\ge 1}I_n\}, \\
B=&\{k\ge 0:(2k+1,2k+2)\subset\cup_{n\ge 1}I_n\}.
\end{align*}
We can choose $a_k,b_k,c_k,d_k$, $k\ge 0$, satisfying the conditions
of Lem\-ma~\ref{lq7} in such a way that
\begin{gather*}
\eta_k=\frac 1{10},\qquad k\in A,\\
\eta_k=e^{-\epsilon(2k+2)(2k+2)},\qquad k\in B,
\end{gather*}
$\Sigma\cap E$ is bounded, and $\Sigma\cap\Lambda(F)=\emptyset$. By
Lemma~\ref{lq7}, we obtain an entire function $G$ of exponential type $2\pi$
such that
\begin{gather}
|G'(\lambda)|\ge c>0,\qquad \lambda\in\Sigma\cap (2k+1,2k+2),\,k\in A,
\notag\\
|G'(\lambda)|\ge c\,e^{-\epsilon(2k+2)(2k+2)},\qquad \lambda\in\Sigma\cap (2k+1,2k+2),\,k\in B,\notag\\
|G(\lambda)|\ge c(1+|\lambda|^2)^{-2},\qquad \lambda\in \Lambda(F).
\label{qm}
\end{gather}

Let $H=FG$. Then $H$ is of exponential type $2\pi$, and by
\eqref{q4}, \eqref{q9}, and \eqref{qm},
$$
|H'(\lambda)|\ge c(1+|\lambda|), \qquad \lambda\in \Lambda(H).
$$

The function
$$
\Phi: \lambda\in \Lambda=\Lambda(H)\mapsto \frac{1}{H'(\lambda)}
$$
belongs to $L^2(\mu)$, where
$$
\mu=\sum_{\lambda\in \Lambda}\delta_\lambda.
$$
Repeating the argument used in the proof of the simple direction of de~Branges' Theorem~\ref{thm-branges}, we see that $\Phi$ annihilates
$\mathcal E(2\pi)$. Therefore, $T(\mu)\ge 2\pi$.

We have
$$
\Lambda=\Lambda(F)\cup\bigcup_{k\ge 0}\{\pm a_k,\pm b_k,\pm c_k,\pm d_k\}.
$$
Put
\begin{gather*}
\Lambda^*=\Lambda(F)\cup\bigcup_{k\in A}\{\pm a_k,\pm b_k,\pm c_k,\pm d_k\}
\cup\bigcup_{k\in B}\{\pm a_k,\pm c_k\},\\
\nu=\sum_{\lambda\in \Lambda^*}\delta_\lambda,\\
V(\lambda)=\left\{
\begin{gathered}
(1+|\lambda|)^{-1},\qquad \lambda\in \Lambda^*,\\
+\infty,\qquad \lambda\notin \Lambda^*.
\end{gathered}
\right.
\end{gather*}

Then $V\in L^2(\nu)$. Let $\epsilon>0$. By de Branges'
Theorem~\ref{thm-branges}, if $\E(\pi+\epsilon)$ is not dense in $C_0(V)$,
then there exists a non-zero function $U\in \mathcal K(\pi+\epsilon,V)$ such that
$\Lambda(U)\subset\Lambda^*$, which contradicts to the Levinson theorem
on the existence of the density of zeros for entire functions of Cartwright
class (see~\cite[Lecture~17]{Levin-book} or \cite[III H2]{Koosis1}).
Therefore, by Bakan's theorem, $\E(\pi+\epsilon)$ is dense in $L^2(\nu)$,
$\epsilon>0$. Thus, $T(\nu)\le \pi$.

In a similar way, we verify that $T(\nu)\ge \pi$, $T(\mu)\le 2\pi$.

Finally, we set $\{x_k\}=\Lambda^*$, and choose $y_k$ such that
\[
\{x_k\}\cup\{y_k\}=\Lambda, \quad {\rm and} \quad
|y_k-x_k|\le e^{-\varepsilon(|x_k|)|x_k|}.
\]
This completes the proof of Theorem~\ref{thm-sharpness}.
\hfill $\square$

\section{Non-classical orthogonal spectral functions. \\
Proof of Theorem~\ref{thm.main}}\label{section-4}

We prove Theorem~\ref{thm.main} in two steps: first we perturb the stably
orthogonal measure $\mu_0$ on $\R$, and then we add a symmetric measure
$\mu_{{\rm i}\R}$ supported by ${\rm i}\R$ with fast decaying tails.

\subsection{Perturbation on the real axis}
Let $\mu_\R$ be a symmetric measure supported by $\R$ such that the integral
\[
\int_\la^\infty \d (\mu_\R - \mu_0)
\]
(conditionally) converges and for some $\delta>0$,
\begin{equation}\label{eq-10}
\int_0^\infty e^{\delta\la} \left| \int_\la^\infty \d (\mu_\R - \mu_0) \right| < \infty\,.
\end{equation}

Then $\widehat{\mu_\R}=\widehat{\mu_\R-\mu_0}+\widehat{\mu_0}$, and
by condition~\eqref{eq-10},
the function $\widehat{\mu_\R-\mu_0}$ has
analytic continuation into the strip $\bigl\{|{\rm Im} z|<\delta\bigr\}$.

By the Gelfand--Levitan theorem~\ref{thm-GL}, the function $\Phi [\mu_0]$ satisfies the Gelfand--Levitan condition
(GL--{\bf i}), and, hence, the function $\Phi [\mu_\R]$ satisfies  (GL--{\bf i}).
Again by the Gelfand--Levitan theorem,
$\mu_\R$ is a spectral
measure of a Sturm--Liouville problem~\eqref{eq_SL}--\eqref{eq_BC} on the
interval $[0,a)$, with the potential $q$ of the same class of smoothness on
$[0,a)$ as the potential $q_0$ that corresponds to $\mu_0$.

Next, we check that the spectral measure $\mu_\R$ is orthogonal. This
follows from our Theorem~\ref{thm_main1} combined with the following claim:

\begin{claim}\label{claim-cl}
There exist positive constants $\delta'>0$ and $C>0$ such that, for all
$x\in\R$,
\[
\mu_\R (I_x)\le C\bigl(\mu_0 (2I_x) + e^{-2\delta'|x|}\bigr),
\]
where $I_x=[x-e^{-\delta'|x|},x+e^{-\delta'|x|}]$ and $2I_x$ is the
concentric interval of twice bigger length.
\end{claim}

\medskip\par\noindent{\em Proof:} Suppose that the claim does not hold;
i.e., for each $n\ge 3$ and each $\delta'>0$, there exists $\la_n$ such that
\[
\mu_\R(I_{\la_n})\ge n\bigl(\mu_0(2I_{\la_n})+e^{-2\delta'|\la_n|}\bigr).
\]
Since the measures $\mu_0$ and $\mu_\R$ are locally finite,
$|\la_n|\to\infty$ when $n\to\infty$. Without loss of generality, we assume that
$\la_n\to +\infty$. Let $\psi(\la)=\int_\la^\infty\d (\mu_\R-\mu_0)$. Then
\begin{multline*}
\psi\bigl(\la_n-e^{-\delta'\la_n}\bigr)-
\psi\bigl(\la_n+e^{-\delta'\la_n}\bigr)
= \mu_\R (I_{\la_n}) - \mu_0 (I_{\la_n})\\ \ge
(n-1) \mu_0 (2I_{\la_n}) + n e^{-2\delta' \la_n}\,.
\end{multline*}
Therefore, at least one of the following two conditions must hold: either
\[
\psi\bigl(\la_n+e^{-\delta'\la_n}\bigr)\le-\frac12(n-1)\mu_0(2I_{\la_n})
-\frac12 ne^{-2\delta'\la_n}\,,
\]
or
\[
\psi\bigl(\la_n-e^{-\delta'\la_n}\bigr)\ge\frac12 (n-1)\mu_0(2I_{\la_n})
+\frac12 n e^{-2\delta'\la_n}\,.
\]
We assume, for instance, that the first case occurs, the second case is
quite similar. Then, for
$\la\in[\la_n+e^{-\delta'\la_n},\la_n+2e^{-\delta'\la_n}]$ we have
\[
\psi(\la)\le\psi\bigl(\la_n+e^{-\delta'\la_n}\bigr)+
\mu_0\bigl(2I_{\la_n}\bigr)\le -e^{-2\delta' \la_n}\,,
\]
whence
\[
\int_{\la_n+e^{-\delta'\la_n}}^{\la_n+2e^{-\delta'\la_n}} |\psi(\la)|
e^{\delta\lambda}\,\d\la\ge
\int_{\la_n+e^{-\delta'\la_n}}^{\la_n+2e^{-\delta'\la_n}}
e^{-2\delta'\la_n+\delta\la}\,\d\la \ge 1
\]
provided that $\delta'\le\delta/3$. Clearly, this contradicts to
\eqref{eq-10}. Hence, the claim. \hfill $\Box$

\subsection{Perturbation on the imaginary axis}\label{add2}

Here, we show that {\em if $\E(a)$ is dense in $L^2(\mu_\R)$, where $\mu_\R$
is a stably orthogonal spectral measure supported by $\R$, then it is also
dense in $L^2(\mu)$}. In the case $a=\infty$, a similar question was studied
by Levitan and Meiman~\cite{LM}, and then by Vul~\cite{Vul}. Later, the same
completeness problem appeared again in Gurarii's work~\cite[Theorem~5]{Gur}
on harmonic analysis in weighted Banach algebras of functions on $\R$ with
asymmetric weights. We cannot use their results since we deal with the case
of finite $a$. Instead, we use an idea from
Yuditskii's work~\cite{Yuditskii} pertaining to the density of polynomials.
The following theorem completes the proof of Theorem~\ref{thm.main}
(condition (GL--{\bf ii}) holds for all measures $\mu_{{\rm i}\R}$
we consider here).

\begin{theorem}\label{thm-5.1}
Let $\mu_\R$ be a non-negative measure on $\R$
satisfying estimate \eqref{eq:growth}, and let $\E(a)$ be stably dense in $L^2(\mu_\R)$
for some $a>0$. Let $\mu_{{\rm i}\R}$ be a measure on ${\rm i}\mathbb R$
such that for some $\delta>0$,
\begin{equation}\label{eq-11}
\int_{\mathbb R}e^{\delta \lambda^2}\,\d\mu_{{\rm i}\R}({\rm i}\lambda)<\infty,
\end{equation}
and let $\mu=\mu_\R+\mu_{{\rm i}\R}$. Then $\E(a)$ is dense in $L^2(\mu)$.
\end{theorem}

\par\noindent{\em Proof:} Suppose that $\E(a)$ is not dense in $L^2(\mu)$,
and denote $X=\clos_{L^2(\mu)}\E(a)$. We will need a lemma which is a
version of a classical result of M.~Riesz and Mergelyan pertaining to the weighted polynomial approximation.

\begin{lemma}\label{lemma-5.1} The elements $f\in X$ extend analytically to
$\C$ with the following estimate: for each $\epsilon>0$, there exists a
constant $C_\epsilon$ {\rm(}independent of $f$ and $z${\rm)} such that
\begin{equation}\label{eq-5.1}
|f(z)|\le C_\epsilon e^{\epsilon|z|^2}\|f\|_{L^2(\mu)},\qquad z\in\C\,.
\end{equation}
\end{lemma}

\noindent{\em Proof:}
Let $h\in L^2(\mu)\ominus X$, $\|h\|_{L^2(\mu)}=1$, and let
$$
H(z)=\int \frac{\overline{h(\lambda)}}{\lambda-z}\,\d\mu(\lambda),\qquad z\in\mathbb C\setminus(\mathbb R\cup{\rm i}\mathbb R).
$$
Then for every $f\in\E(a)$ we have
$$
\int \frac{f(\lambda)-f(z)}{\lambda-z}\overline{h(\lambda)}\,\d\mu(\lambda)=0,\qquad z\in\mathbb C,
$$
and, hence,
\begin{equation}\label{fcd1}
f(z)=\frac 1{H(z)}\int
\frac{f(\lambda)\overline{h(\lambda)}}{\lambda-z}\,
\d\mu(\lambda),\qquad \lambda\in\mathbb C\setminus(\mathbb R\cup {\rm i}\mathbb R).
\end{equation}

Since the function $z\mapsto \int
\frac{f(\lambda)\overline{h(\lambda)}}{\lambda-z}\,
\d\mu(\lambda)$ is bounded and the function $H$ is of at most linear growth
in $\overline{\Omega}$, where
$\Omega=\{z\in\mathbb C:\dist(z,\mathbb R\cup {\rm i}\mathbb R)>1\}$,
we obtain that $H$ and, hence, $f$ are in the Nevanlinna class (see
\cite[Section II.5]{GAR}) in $\Omega$.
The function $f$ is of exponential type in the plane, and we conclude that
$$
\log |f(z)|\le\int_{\partial\Omega}\log|f(\lambda)|\,
\omega(z,\d \lambda,\Omega),
$$
where $\omega(z,E,\Omega)$ is the harmonic measure of $E\subset\partial \Omega$ in $\Omega$ with respect to
$z\in\Omega$.
By \eqref{fcd1}, $|f|\le \|f\|_{L^2(\mu)}/|H|$ on $\partial\Omega$.
Since
$\log|H|\in L^1(\omega(z,\d \lambda,\Omega))$, $z\in\Omega$,
for every $\epsilon>0$, we can find $c_{\epsilon}>0$ such that
$$
|f(z)|\le c_{\epsilon} \|f\|_{L^2(\mu)}e^{\epsilon |z|^2}, \qquad z\in L\cap\Omega_1,
$$
where $L=\{re^{i\theta}:0<r<\infty,\, |\sin 2\theta|> 1/2\}$,
$\Omega_1=\{z\in\mathbb C:\dist(z,\mathbb R\cup {\rm i}\mathbb R)\ge 2\}$.

Once again, since $f\in\E(a)$, the Phragm\'en-Lindel\"of principle applied to $f$
in each of the four sectors of $\mathbb C\setminus L$
gives
$$
|f(z)|\le c_\epsilon\|f\|_{L^2(\mu)} e^{\epsilon |z|^2}, \qquad z\in \mathbb C.
$$
Passing to the limit, we obtain the same inequality for all $f\in X$.
\mbox{} \hfill $\Box$

\medskip We resume the proof of Theorem~\ref{thm-5.1}, and define a linear
operator $K$ on $X$ by the relation
\[
\langle Kf,g \rangle_X  = \int_{{\rm i}\R} f\bar g\,\d\mu_{{\rm i} \R}.
\]
Since the form on the right-hand side is bounded, we see that the operator
$K$ is bounded. Furthermore, since
$$
\langle (I-K)f,f \rangle_X = \int_\R |f|^2\,\d\mu_\R\,,
$$
we see that $0\le K \le I$ in the operator sense.

\begin{lemma}\label{Lemma-5.2} The operator $K$ is compact.
\end{lemma}

\noindent{\em Proof:} Let $f_n\in X$ tend weakly to $0$. By the
Banach--Steinhaus uniform boundedness principle,
$\sup_n\|f_n\|_{L^2(\mu)}<\infty$. Therefore, by~\eqref{eq-5.1},
the family of entire functions $\{f_n\}$ is equicontinuous on every compact subset of
$\C$. Again applying~\eqref{eq-5.1}, we see that $f_n$ tend to $0$ pointwise
on $\C$, and hence, uniformly on compact subsets of $\C$. Furthermore,
\[
\left| \langle K f_n, g \rangle_X \right|^2 \le
\| g \|^2_{L^2(\mu)}\cdot
\int_{\R} |f_n({\rm i}y)|^2\, \d\mu_{{\rm i}\R} ({\rm i}y)\,.
\]
By the dominated convergence theorem (it can be applied due to
estimates~\eqref{eq-5.1} and \eqref{eq-11}), the integral in the right-hand side
tends to $0$ when $n\to\infty$, whence,
$$
\lim_{n\to\infty}\ \sup_{\|g\|_{L^2(\mu)}\le 1}\ |\langle Kf_n,g\rangle_X|=0\,,
$$
proving the compactness of the operator $K$. \hfill $\Box$

\medskip We proceed with the proof of Theorem~\ref{thm-5.1}, and denote by
$\sigma(K)$ the spectrum of the operator $K$. First, suppose that
$1\notin \sigma (K)$. Then the operator $I-K$ is invertible. Therefore, for each $g\in\E(a)\subset X$, we have
\begin{multline*}
\|g\|^2_{L^2(\mu)}=\|g\|^2_X=\langle (I-K)^{-1}(I-K)g,g \rangle_X \\
\le \|(I-K)^{-1}\|\cdot\langle (I-K)g, g \rangle_X =
\|(I-K)^{-1}\|\cdot\|g\|^2_{L^2(\mu_\R)}\,,
\end{multline*}
and by~\eqref{eq-5.1},
\[
|g(z)|\le C_\epsilon e^{\epsilon|z|^2}\,\|(I-K)^{-1}\|^{1/2}\cdot
\|g\|_{L^2(\mu_\R)}\,,\qquad z\in\C\,.
\]
In particular, $\kappa=\inf\bigl\{\|g\|_{L^2(\mu)}:|g({\rm i})|=1\bigr\}>0$.
However, this contradicts to the stable density of $\E(a)$ in $L^2(\mu_\R)$.
Indeed, for any function $f\in\E(a)$ with $f({\rm i})=1$, and any
$h\in\E(a)$, we have
\begin{multline*}
\bigl\|h(x)-f(x)(x-{\rm i})^{-1}\bigr\|_{L^2 (\mu_2)}\asymp
\|h(x)(x-{\rm i})-f(x)\|_{L^2(\mu)}\\
\ge\inf\{\|g\|_{L^2(\mu)}:|g({\rm i})|=1\}=\kappa>0\,,
\end{multline*}
where, as above, $\d\mu_2(\la)=(1+|\la|)^2\,\d\mu(\la)$. Recalling
Theorem~\ref{la}, we conclude that $\E(a)$ is not dense in $L^2(\mu_2)$, and
hence, is not stably dense in $L^2(\mu)$.

\medskip Now, we suppose that $1\in\sigma (K)$. In this case, the operator
$I-K$ is not invertible, but its kernel is finite dimensional,
and
we can modify the previous argument.

Denote $X_0=\operatorname{ker}(I-K)$,
$X_1=X\ominus X_0$, $N=\operatorname{dim}(X_0)$, and $K_1=K|\,X_1$. Note
that $\|K_1\|<1$. Take a basis $g_1,\ldots, g_N$ in $X_0$, and choose
points $x_1,\ldots,x_N$ on $\R$ such that the matrix
$\bigl[g_j(x_k)\bigr]_{1\le j,k\le N}$ is non-degenerate (this choice is
possible due to the linear independence of the functions $g_1,\ldots,g_N$).
Denote $\nu=\sum_k\delta_{x_k}$. Then $|g(i)|\le c\|g\|_{L^2(\nu)}$,
$g\in X_0$.

Next, set $\widetilde\mu_\R=\mu_\R+\nu$. Let $f=f_0+f_1\in X$, $f_0\in X_0$,
$f_1\in X_1$, $\|f\|_{L^2(\widetilde\mu_\R)}\le 1$. Since
$(I-K)(f_0+f_1)=(I-K_1)f_1$, and
$$
\langle (I-K)(f_0+f_1),f_0+f_1\rangle_X=\langle(I-K_1)f_1,f_1\rangle_X,
$$
we obtain, as above, that
$\|f_1\|_{L^2(\mu)}\le c$. By Lemma~\ref{lemma-5.1},
$\|f_1\|_{L^2(\nu)}\le c_1$ and $|f_1(i)|\le c_2$. Therefore,
$\|f_0\|_{L^2(\nu)}\le 1+c_1$, $|f_0(i)|\le c_3$, and $|f(i)|\le c_2+c_3$.

As a result, we obtain that $|f(i)|\le c\|f\|_{L^2(\widetilde\mu_\R)}$ for
each $f\in \E (a)$. By Lemma~\ref{lemma-stable} and Theorem~\ref{la}, this
contradicts to the stable density of $\E (a)$ in $L^2(\mu_\R)$.
\mbox{} \hfill $\Box$

\appendix
\section{Cases when the type is explicitly computable}\label{appendix-1}

There are several cases when the type $T(\mu)$ can be explicitly computed,
or at least estimated. Here, we briefly list some of these cases.

\subsection{Measures of zero type}\label{xx2}

If the tails of the measure $\mu$ decay sufficiently rapidly or if there are
large gaps in the support of $\mu$, then the measure has zero type. A useful
sufficient condition that deals with these two kinds of behavior is due to
de~Branges~\cite[Theorem~63]{deBr}; an equivalent result was obtained by
Beurling, see \cite[Section~VII.A.2]{Koosis1}.
{\em Let $K\colon \R\to [1, +\infty]$, $\log K$ be uniformly continuous,
and let
\[
\int_\R \frac{\log K(t)}{t^2+1}\,\d t = \infty\,.
\]
If $\displaystyle \int_\R K\,\d\mu<\infty$, then $T(\mu)=0$}.

\subsection{Measures of infinite type}\label{xx3}

A useful sufficient condition for $T(\mu)=\infty$ is due to Krein (and goes back to Szeg\"o).
{\em Suppose that $\mu$ has a bounded density $\mu'$ with respect to Lebesgue measure, and that
\[
\int_\R \frac{\log \mu'(t)}{t^2+1}\, \d t > -\infty\,.
\]
Then $T(\mu)=\infty$}.

Another useful result follows from a theorem of
Duffin and Schaeffer~\cite{DuSc}: if the measure $\mu$ is relatively dense
with respect to Lebesgue measure, then it must have a positive type. More
precisely, {\em suppose that for some $L<\infty$ and $\delta>0$,
$$
\mu [x-L, x+L]\ge \delta,\qquad x\in\mathbb R.
$$
Then $T(\mu) \ge \frac{2\pi}L$.} This can be regarded as a certain
stability of the infinite type of Lebesgue measure.

\subsection{Measures of positive type supported by discrete separated sets}\label{appendix-1.3}

If the measure $\mu$ is supported by the set of the integers $\Z$, then the Fourier
transforms of the measures $g\,\d\mu$ are $2\pi$-periodic functions, whence
$T(\mu)\le \pi$. A result of Koosis~\cite{Koosis2} yields that {\em if
$\mu=\sum_\Z\omega(n)\delta_n$, where $\delta_n$ is the point mass at $n$,
and $\omega\colon\Z\to\R_+$ is an arbitrary sequence such that
$\sum_{n\in\Z}\frac{\omega(n)}{1+n^2}<\infty$ and
$\sum_{n\in\Z}\frac{\log\omega(n)}{1+n^2}>-\infty$, then $T(\mu)=\pi$}. This
is a deep result which readily yields one of the equivalent forms of the
Beurling--Malliavin multiplier theorem.  It is worth mentioning that this
fact has a much simpler proof in the case when $\omega$ satisfies additionally
some regularity assumptions, for instance, if it is an even
non-decreasing sequence.

Another deep result is a recent theorem of Mitkovski and
Polto\-ratski~\cite{PolMit} which, in its turn, goes back to de~Branges. It
yields that if $\Lambda=\{\la\}\subset\R$ is a separated sequence of points
and $\mu=\sum_\Lambda\delta_{\la_n}$, then
$T(\mu)=\pi\mathcal D_*(\Lambda)$, where $\mathcal D_*(\Lambda)$ is the
{\em lower Beurling--Malliavin density} of the sequence $\Lambda$.

\section{Stable and unstable density}\label{appendix-2}

Here, we will discuss measures $\mu$ such that $\E(a)$ is dense but not
stably dense in $L^2(\mu)$. Our discussion is close to the one in~\cite[Appendix~1]{BoSo}
where we dealt with the weighted polynomial approximation.

As above, given a real $t$, we define the measure $\mu_t$ and the weight $W_t$ as follows:
\[
\d\mu_t(\la)=(1+|\la|)^t\,\d\mu(\la)\,, \quad W_t(\la) = W(\la) (1+|\la|)^{-t}\,.
\]

\begin{lemma}\label{cor-singular}
Let $a>0$, and let $\mu$ be a non-negative measure on $\R$ of at most
polynomial growth such that $\supp\mu$ does not contain the zero set of any
function from the Krein class $\mathcal K(a)$. Then $\E(a)$ is stably dense
in $L^2(\mu)$.
\end{lemma}

\par\noindent{\em Proof:} Suppose that for some $t\in\R$, $\E(a)$ is not
dense in $L^2(\mu_t)$. For some $N<\infty$, the weight function $W$
defined by
$$
W(x)=\begin{cases}
(1+|x|)^{-N},\qquad x\in\supp\mu,\\
+\infty,\qquad x\not\in\supp\mu,
\end{cases}
$$
belongs to $L^2(\mu_t)$. By Bakan's
Theorem~\ref{lb}, $\E(a)$ is not dense in $C_0(W)$. By
de Branges' Theorem~\ref{thm-branges}, there exists $f\in \mathcal K(a)$ such
that \[ \Lambda(f)\subset\{x:W(x)\not=+\infty\}=\supp\mu\,.\]
\hfill $\Box$

\begin{definition}\label{def-a_singular}
Let $a>0$. We say that the measure $\mu$ is $a$-singular, if there exist
real $t$ and $s$, $t<s$, such that $\E(a)$ is dense in $L^2(\mu_t)$ and is
not dense in $L^2(\mu_s)$. Similarly, we say that the weight $W$ is $a$-singular,
if there exist real $t$ and $s$, $t<s$, such that $\E (a)$ is dense in $C_0(W_t)$ and
is not dense in $C_0(W_s)$.
\end{definition}

\begin{lemma}\label{lemma-singular}\mbox{}

\smallskip\par\noindent {\rm (I)}
Suppose $W$ is an $a$-singular weight. Then there exists
a function $B$ of Krein's class $\mathcal K(a)$ such that $\{\la\colon W(\la)\ne \infty\} = \Lambda (B)$.

\smallskip\par\noindent {\rm (II)}
Suppose $\mu$ is an $a$-singular measure. Then there exists
a function $B$ of Krein's class $\mathcal K(a)$ such that $\operatorname{supp}(\mu)=\Lambda (B)$.
\end{lemma}

\par\noindent{\em Proof of Lemma~\ref{lemma-singular}:} The proof will be similar to
the previous one.

\smallskip\par\noindent {\rm (I)} Let $W$ be an $a$-singular weight. Since the weights $W$ and $W_t$
are finite on the same set of points, we may assume that $\E (a)$ is dense in $C_0(W)$
and is not dense in $C_0(W_1)$. Then by de Branges' Theorem~\ref{thm-branges}, there is a function
$B$ in $\mathcal K(a)$ such that
\[
\sum_{\la\in\Lambda (B)} \frac{W(\la)}{(1+|\la|)|B'(\la)|} < \infty\,.
\]
Clearly, $\{\la\colon W(\la)\ne \infty\} \supset \Lambda (B)$. Suppose that
$\{\la\colon W(\la)\ne \infty\} \supsetneqq \Lambda (B)$, take a point
$\lambda_0\in \{\la\colon W(\la)\ne \infty\} \setminus \Lambda (B)$, and consider
the entire function $B_1(z)=(z-\la_0)B(z)$. Clearly, this is again a function of
Krein's class $\mathcal K(a)$, $\Lambda(B_1) = \Lambda(B) \cup \{\la_0\}$,
and for $\la\in\Lambda(B)$, we have $|B_1'(\la)| = |\la-\la_0| |B'(\la)|$. Therefore,
\[
\sum_{\la\in\Lambda (B_1)} \frac{W(\la)}{|B_1'(\la)|} < \infty\,,
\]
and by the other half of de Branges' Theorem, $\E (a)$ is dense in $C_0(W)$,
which contradicts our assumption.

\smallskip\par\noindent {\rm (II)} Let $\mu$ be an $a$-singular measure. As above, we assume
that $\E (a)$ is dense in $L^2(\mu)$ and is not dense in $L^2(\mu_1)$. By Bakan's theorem~\ref{lb},
there exists an $a$-singular weight $ W\in L^2(\mu) $. Therefore, applying the first part of the
lemma, we see that the support $\Lambda$ of the measure $\mu$ is contained in the (discrete) set $\{\la\colon W(\la)\ne \infty\}$ which coincide with
the zero set of a function of
Krein's class. In particular,
\[
\sum_{\la\in\Lambda} \frac1{(1+|\la|)^2} < \infty\,.
\]

Introduce an auxiliary weight $V$,
\[
V(\la) =
\begin{cases}
\mu \{\la\}^{-1/2}, & \la\in\Lambda, \\
\infty, & {\rm otherwise}\,.
\end{cases}
\]
Then for every function $\phi$ with a compact support,
\begin{gather*}
\|\phi\|^2_{C_0(V)}=\max_{\la\in\Lambda}|\phi(\la)|^2 \mu\{\la\}\le
\sum_{\la\in\Lambda}|\phi(\la)|^2\mu\{\la\}=\|\phi\|^2_{L^2(\mu)},\\
\|\phi\|^2_{L^2(\mu_1)}\le\max_{\la\in\Lambda}|\phi(\la)|^2(1+|\la|)^3 \mu\{\la\}\cdot\sum_{\la\in\Lambda}\frac1{(1+|\la|)^2}\le
C(\Lambda) \|\phi\|^2_{C_0(V_{3/2})}\,.
\end{gather*}
Thus, the  weight $V$ is $a$-singular, and the first part of the lemma
completes the proof. \hfill $\Box$

\medskip
We also use the following lemma:

\begin{lemma}\label{lemma-stable}
Suppose that $\E(a)$ is stably dense in $L^2(\mu)$,
$\la_k\in\R$, $1\le k\le N$, and let
$\displaystyle\widetilde\mu=\mu+\sum\limits_{1\le k\le N}\delta_{\la_k}$.
Then $\E(a)$ is dense in $L^2(\widetilde\mu)$.
\end{lemma}

\par\noindent{\em Proof:}
By Bakan's Theorem, there exists a weight $W\in L^2(\mu)$ such that $\E(a)$
is stably dense in $C_0(W)$. Let $\widetilde W\le W$ be a (lower
semi-continuous) weight such that $W(\la_k)<+\infty$, $1\le k\le N$,
$\widetilde W=W$ outside of a finite interval $I$.

We claim that $\E(a)$ is dense in $C_0(\widetilde W)$. Otherwise, by
de Branges' Theorem~\ref{thm-branges}, we find a Krein class function
$\widetilde B$ such that
\[
\sum_{\la\in\Lambda(\widetilde B)} \frac{\widetilde W (\la)}{| \widetilde B'(\la)|} < \infty\,.
\]
Set $\Lambda=\Lambda(\widetilde B)\cap I$, and let
\[
B(z)\stackrel{\rm def}=\frac{\widetilde B(z)}
{\prod\limits_{\la\in\Lambda}(z-\la)}\,.
\]
Clearly, $B$ is a Krein class function, and for some $s<\infty$, we
have $|B'(\la)|\asymp(1+|\la|)^{-s}|\widetilde B'(\la)|$ for
$\la\in\Lambda (B)$. Thereby,
\[
\sum_{\la\in\Lambda(B)} \frac{W(\la)}{(1+|\la|)^s|B'(\la)|} < \infty\,,
\]
which by de Branges' theorem contradicts to the stable density of $\E(a)$ in
$C_0(W)$.

Thus, $\E(a)$ is dense in $C_0(\widetilde W)$. Since
$\widetilde W\in L^2(\widetilde\mu)$, using again Bakan's Theorem, we
conclude that $\E(a)$ is dense in $L^2(\widetilde \mu)$, completing the
proof. \hfill $\Box$

\medskip It is not difficult to prove the inverse statement:
{\em if for any points $\la_1,\ldots,\la_N\in\R$, $\E(a)$ is dense in
$L^2(\widetilde\mu)$, where $\widetilde\mu=\mu+\sum_{k}\delta_{\la_k}$,
then  $\E(a)$ is stably dense in $L^2(\mu)$.}

\section{Nazarov's construction of spectral measures}\label{appendix-5}

We describe an elegant Nazarov's construction of a wide class of spectral
measures supported by $\R$. This construction is based on a distorted
Poisson formula.

\begin{definition}\label{c11}
Denote by $\Gamma$ the set of $C^\infty$-diffeomorphisms $X$ of $\R$
satisfying the following two conditions:

\smallskip\par\noindent {\rm (I)} $X'(t)\to 1$ as $|t|\to\infty$;

\smallskip\par\noindent {\rm (II)} There exists a sequence of {\rm(}strictly{\rm)} positive
numbers $s_k$, $\displaystyle\lim_{k\to\infty}s_k=\infty$, such that for every
$k\ge 2$, $X^{(k)}(t)=O(|t|^{-s_k})$ as $|t|\to\infty$.
\end{definition}

\medskip Given a diffeomorphism $X\in\Gamma$ and given $c>0$, we define the
measure
\begin{equation}\label{eq-Nazarov}
\mu=\mu_X=\sum_{k\in\Z}X'(ck)\delta_{X(ck)} = \sum_{\la\colon Y(\la)\in c\Z} \frac{\delta_\la}{Y'(\la)} \,.
\end{equation}
Here and below, $Y=X^{-1}$ is the inverse diffeomorphism. By $\widehat\mu$
we denote the distributional Fourier transform of $\mu$, i.e.,
$\langle\wh\mu,\phi\rangle=\langle\mu,\wh\phi\rangle $; as above,
\[
\wh \phi (\la) = \int_\R \phi (x) e^{-{\rm i} \la x}\,\d\la\,.
\]

\begin{theorem}\label{lemma-Nazarov}
There exists a function $M\in C^\infty (-2\pi c^{-1}, 2\pi c^{-1})$ such that
$\widehat\mu=\delta_0+M$ on $(-2\pi c^{-1},2\pi c^{-1})$.
\end{theorem}

\medskip
We start with an obvious claim:

\begin{claim}\label{claim-1}
Let $s_k$, $k\ge k_0$, be any sequence of {\rm(}strictly{\rm)} positive numbers satisfying
$\displaystyle\lim_{k\to\infty}s_k=\infty$. Then the sequence
$$
S_m = \inf\Bigl\{\sum_j s_{k_j}\colon  \sum_j k_j\ge m\Bigr\}
$$
also satisfies $\displaystyle \lim_{m\to\infty}S_m = \infty$.
\end{claim}

Next, we prove

\begin{claim}\label{claim-2}
If $X\in\Gamma$, then $Y=X^{-1}\in\Gamma$.
\end{claim}

\par\noindent{\em Proof:}
Condition (I) implies that $X(t)\asymp t$ for large $t$
and, thereby, $Y(t)\asymp t$ for large $t$ as well. Then condition (I) for $Y$ follows immediately from the identity $Y'=\frac 1{X'\circ Y}$.
Differentiating this identity $k-1$ times, we conclude that $Y^{(k)}$ is a
finite linear combination of terms of the kind
\[
(X'\circ Y)^{-m}\prod_j X^{(k_j)}\circ Y, \quad  {\rm with\ } m,k_j\ge 1, \quad \sum_j(k_j-1)= k-1.
\]
Now applying Claim~\ref{claim-1} to the sequence $\{s_{k+1}\}_{k\ge 1}$, we
conclude that each such term is $O(|t|^{-S_k})$ where $\{S_k\}_{k\ge 2}$ is a
sequence of positive numbers tending to infinity. \hfill $\Box$

\medskip\par\noindent{\em Proof of Theorem~\ref{lemma-Nazarov}:} It suffices
to deal only with the case $c=1$ (otherwise, just replace $X$ by
$t\mapsto c^{-1}X(c\cdot t)$). We fix any $0<a<b<2\pi$, and take
a function $\phi$ in the Schwartz class $\mathcal S$ such that
$\wh\phi=1$ on $(-a,a)$ and $\wh\phi=0$ outside $(-b,b)$. Clearly,
$\wh\mu=\wh\mu\wh\phi=\wh{\mu*\phi}$ on $(-a,a)$. Taking into account
that the point mass at the origin is the Fourier transform of the constant
$\tfrac{1}{2\pi}$, it remains to prove that
$(\mu*\phi)(t)-\tfrac1{2\pi}$ decays faster than any power of $|t|$ as
$|t|\to\infty$.

To this end, fix large $t\in\R$ and consider the function
$\psi(s)=\phi(t-X(s))X'(s)$. Since all the derivatives of $X$ are bounded and
$X(s)\asymp s$ for large $s$, the function $\psi$ belongs to $\mathcal S$.
Moreover, we have $\|\psi^{(k)}\|_{L^1}\le C(k,X,\phi)$ independently of
$t$, which allows us to conclude that
$|\widehat\psi(\ell)|\le C(N,X,\phi)(1+|\ell|)^{-N}$ for any $N>0$, independently
of $t$ and $\ell$.

Next,
$$
(\mu*\phi)(t)=\sum_{k\in\Z}X'(k)\phi(t-X(k))=\sum_{k\in\Z}\psi(k)=
\frac1{2\pi} \sum_{k\in\Z}\wh\psi \left( 2\pi k \right)\,.
$$
Since
$$
\wh\psi(0)=\int_\R\psi=\int_\R \phi(t-X(s))X'(s)\,\d s=\int_\R\phi=\wh\phi(0)=1\,,
$$
we need only to prove that
$\sum\limits_{\ell\in\Z\setminus\{0\}}|\wh\psi(2\pi \ell)|$ decays
faster than any power of $|t|$ as $|t|\to\infty$.

Note that if $X$ is the identical map, this property holds because the support of
$\wh\psi$ is contained in $(-b,b)$ with $b<2\pi$. We need to show that this
property is preserved under distortion of $\phi$ by any diffeomorphism $X$
which is asymptotically close to the identity.

We put $Y=X^{-1}$, $R(t)=Y(t)-t$, and write
\begin{multline*}
\wh\psi(2\pi\ell)=\int_\R \phi(t-X(s))X'(s)e^{-2\pi {\rm i}\ell s}\,\d s \\
=\int_\R\phi(t-s)e^{-2\pi {\rm i}\ell Y(s)}\,\d s
=\int_\R\phi(s)e^{-2\pi {\rm i}\ell Y(t-s)}\,\d s \\
=e^{-2\pi {\rm i}\ell R(t)-2\pi {\rm i}\ell t}
\int_\R\phi(s)e^{2\pi {\rm i}\ell s+2\pi {\rm i}\ell [R(t)-R(t-s)]}\,\d s\,.
\end{multline*}

\begin{claim}\label{claim:dop} There exists $0<\delta<\tfrac14$ such that for large $N$,
\[
\int_\R\phi(s)e^{2\pi {\rm i}\ell s+2\pi {\rm i}\ell [R(t)-R(t-s)]}\,\d s =
O\bigl( |t|^{-N} \bigr), \qquad |t|\to\infty,
\]
uniformly in $1\le \ell \le |t|^\delta$.
\end{claim}
Since $\sum\limits_{\ell\colon|\ell|>|t|^\delta}|\wh\psi(2\pi\ell)|$
decays faster than any power of $|t|$ as $|t|\to\infty$, proving the claim,
we prove the theorem.

\medskip\par\noindent{\em Proof of Claim~\ref{claim:dop}:}
Fix $N\in(1,\infty)$ and write the Taylor formula:
\begin{multline*}
Q(t,s) \stackrel{\rm def}= R(t)-sR'(t)-R(t-s)
\\
=\sum_{k=2}^{p-1}a_k(t)s^k+O(a_p(t)s^p),\qquad |t|\to\infty,\,
|s|<|t|/2.
\end{multline*}
By Claim~\ref{claim-2} and the definition of the class $\Gamma$,
we can find $0<\delta<\tfrac14$, $p<\infty$, and $C<\infty$
such that for large $|t|$,
\begin{align}
\label{eq:dop5}
|a_k(t)|&\le C|t|^{-4\delta},\qquad 2\le k\le p-1,\\
\nonumber
|a_p(t)|&\le C|t|^{-2N}.
\end{align}
We fix $0<\epsilon<\frac{\delta}p$, and set
$$
A(t,s)=\sum_{k=2}^{p-1}a_k(t)s^k.
$$
Then for $|s|<|t|^\epsilon$, $|t|>1$, we have
\begin{align}\label{eq:dop6}
|Q(t,s)-A(t,s)|&\le c|t|^{-2N+\delta}, \\
\label{eq:dop7}
|A(t,s)|&\le c|t|^{-3\delta}.
\end{align}

Furthermore,
$$
\int_{\R}\phi(s)e^{-2\pi {\rm i}su}\,\d s=0, \qquad |u|\ge b\,,
$$
and hence for large $|t|$ we have
\[
\int_{\R}s^k\phi(s)e^{2\pi {\rm i}s\ell(1+R'(t))}\,\d s =
0, \qquad \ell\in\mathbb Z\setminus\{0\},\, k\ge 0\,.
\]
Since $\phi\in\mathcal S$, we have
\begin{equation}\label{eq:dop3}
|\phi(t)|=O(t^{-1-3N/\epsilon}),\qquad |t|\to \infty,
\end{equation}
whence,
\begin{equation}\label{eq:dop4}
\int_{|s|<|t|^\epsilon}s^k\phi(s)e^{2\pi {\rm i}s\ell(1+R'(t))}\,\d s
= O(t^{-2N}),\qquad |t|\to\infty\,,
\end{equation}
uniformly in $0\le k\le pN/\delta$.

Using~\eqref{eq:dop3} once more, we get
\begin{multline*}
\int_\R\phi(s)e^{2\pi {\rm i}\ell s+2\pi {\rm i}\ell [R(t)-R(t-s)]}
\,\d s
\\
\!\!\!\!\!\!\!\!\!\!= \int_{|s|<|t|^\epsilon} \phi(s)e^{2\pi {\rm i}\ell s+2\pi {\rm i}\ell [R(t)-R(t-s)]}
\,\d s + O\bigl( |t|^{-N} \bigr) \\
= \int_{|s|<|t|^\epsilon} \phi(s)e^{2\pi {\rm i}\ell s (1+R'(t))}
e^{2\pi {\rm i}\ell Q(t,s)}
\,\d s + O\bigl( |t|^{-N} \bigr), \qquad |t|\to\infty.
\end{multline*}

Next, choosing an integer $q$, $\frac{N+1}{2\delta}<q<
\frac{N}{\delta}$, we expand
\begin{gather*}
e^{2\pi {\rm i} \ell Q}
= \sum_{j=0}^{q-1} \frac{\left( 2\pi {\rm i} \ell Q \right)^j}{j!}
+ O\bigl( |\ell Q|^q \bigr)
=
\sum_{j=0}^{q-1} \frac{\left( 2\pi {\rm i} \ell \left[ A + (Q-A) \right]
\right)^j}{j!} + O\bigl( |\ell Q|^q \bigr)
\\
=
\sum_{j=0}^{q-1} \frac{\left( 2\pi {\rm i} \ell \right)^j}{j!} A^j  +
\sum_{j=1}^{q-1} O\bigl( |\ell|^j \sum_{m=0}^{j-1} |A|^m |Q-A|^{j-m} \bigr)
+ O\bigl( |\ell Q|^q \bigr)
\end{gather*}
for $|t|\to \infty$, uniformly in $|s|\le |t|^\epsilon$ and $|\ell|\le |t|^\delta$.
Using estimates~\eqref{eq:dop6}-\eqref{eq:dop7}, and recalling that $N>4\delta$,
and that $\tfrac12 (N+1) < \delta q < N$, we get
\[
e^{2\pi {\rm i} \ell Q} = \sum_{j=0}^{q-1} \frac{\left( 2\pi {\rm i} \ell \right)^j}{j!} A^j
+ O\bigl( |t|^{-N-2\delta}\bigr), \qquad |t|\to\infty\,,
\]
again, uniformly in $|s|\le |t|^\epsilon$ and $|\ell|\le |t|^\delta$.
Then we substitute this expansion into the integral we are estimating:
\begin{multline*}
\int_\R\phi(s)e^{2\pi {\rm i}\ell s+2\pi {\rm i}\ell [R(t)-R(t-s)]}
\,\d s
\\
= \sum_{j=0}^{q-1} \frac{\left( 2\pi {\rm i} \ell \right)^j}{j!}
\int_{|s|<|t|^\epsilon} \phi(s) A(t,s)^j e^{2\pi {\rm i}\ell s (1+R'(t))}\,\d s
+ O\bigl( |t|^{-N} \bigr)
\\
\stackrel{\eqref{eq:dop4}}=
\sum_{j=1}^{q-1} \frac{\left( 2\pi {\rm i} \ell \right)^j}{j!}
\int_{|s|<|t|^\epsilon} \phi(s) A(t,s)^j e^{2\pi {\rm i}\ell s (1+R'(t))}\,\d s
+ O\bigl( |t|^{-N} \bigr)
\end{multline*}
for $|t|\to\infty$, uniformly in $1\le |\ell|\le|t|^\delta$. At last, note that
the functions $s\mapsto A(t, s)^j$ are polynomials of degree at most $(p-1)(q-1)$ which
is less than $pN/\delta$ by the choice of $q$. Hence, integrating these polynomials
we can apply estimates~\eqref{eq:dop4}. By~\eqref{eq:dop5}, the coefficients of these polynomials
are bounded by $O(|t|^{-4\delta j})$. Hence,
\begin{multline*}
\sum_{j=1}^{q-1} \frac{\left( 2\pi {\rm i} \ell \right)^j}{j!}
\int_{|s|<|t|^\epsilon} \phi(s) A(t,s)^j e^{2\pi {\rm i}\ell s (1+R'(t))}\,\d s
\\
= O\bigl( \sum_{j=1}^{q-1} |\ell|^j |t|^{-4\delta j}\, |t|^{-2N} \bigr)
= O(|t|^{-N}), \qquad |t|\to\infty\,,
\end{multline*}
uniformly in $1\le |\ell| \le |t|^\delta$.
This proves the claim and finishes off the proof of Theorem~\ref{lemma-Nazarov}.
\hfill $\Box$

\section{Stably orthogonal spectral measures}\label{appendix-6}

In this appendix, we show that Nazarov's construction of spectral measures
described in Appendix~\ref{appendix-5} gives a wide class of explicitly
defined discrete stably orthogonal spectral measures.
Let us take a diffeomorphism $X$ satisfying Definition~\ref{c11}. Suppose that $X$ is an
odd function. Then the measure $\mu$ defined in~\eqref{eq-Nazarov} is symmetric
with respect to the origin, and hence, is a spectral measure of a
Sturm--Liouville problem with a potential $q\in C^\infty[0,\pi c^{-1})$.
Now, we show that if
\begin{equation}\label{eq-Nazarov2}
\lim_{t\to+\infty} \bigl( X(t)-t \bigr) = +\infty\,,
\end{equation}
then the spectral measure $\mu$ is stably orthogonal.

For a discrete set $\Lambda\subset \R$, we denote
\[
n_\Lambda(t)=\#\{\Lambda\cap[-t,t]\},\qquad
N_\Lambda(R)=\int_1^R \frac{n_\Lambda (t)}{t}\,\d t\,.
\]
Let $\Lambda=X(c\Z)$. Then by condition~\eqref{eq-Nazarov2},
\begin{equation}\label{eq-N}
\lim_{t\to \infty} \left( n_\Lambda (t) - 2t c^{-1} \right) = -\infty\,.
\end{equation}
By Lemma~\ref{cor-singular}, the following lemma does the job.

\begin{lemma}\label{lemma-Krein_zeros}
Suppose that the set $\Lambda$ satisfies condition~\eqref{eq-N}. Then
$\Lambda$ cannot contain the zero set of any function $F\in \mathcal K(\pi c^{-1})$.
\end{lemma}

\par\noindent{\em Proof:} Suppose that for some function $F\in \mathcal K(\pi c^{-1})$, we
have $\Lambda(F)\subset\Lambda$, where $\Lambda(F)$ is the zero set of $F$.
Take a polynomial $P$ with real coefficients and with simple zeros in
$\R\setminus\Lambda(F)$. Then the entire function $G=P\cdot F$ again belongs
to the Krein class $\mathcal K(\pi c^{-1})$, and $|G'(\la)|\ge c(1+|\la|^2)$, provided
that the degree of $P$ is big enough.

Take any $b<\pi c^{-1}$, and write the Lagrange interpolation formula,
\[
\frac{\sin (bz)}{G(z)} = \sum_{\Lambda (G)} \frac{\sin (b\la)}{G'(\la)(z-\la)}\,.
\]
Then
\[
|\sin(bz)|\le\frac{|G(z)|}{|\Im z|}\sum_{\Lambda(G)}\frac1{1+|\la|^2}=
C\frac{|G(z)|}{|\Im z|}\,.
\]
Since the constant $C$ does not depend on the choice of $b<\pi c^{-1}$, letting
$b\to\pi c^{-1}$, we obtain $|\sin(\pi c^{-1}z)|\le C|\Im z|^{-1}\,|G(z)|$, and then
\begin{multline*}
\log |\sin (\pi c^{-1}z)| \le \log |G(z)| + \log \frac1{|\Im z|} + C' \\
= \log |F(z)| + \log |P(z)| + \log \frac1{|\Im z|} + C'\,.
\end{multline*}
Letting $z=Re^{{\rm i}\theta}$, integrating over $\theta$, and using
Jensen's formula, we get
\begin{multline*}
2c^{-1} R = \frac1{2\pi} \int_0^{2\pi} \log |\sin (\pi c^{-1}Re^{{\rm i}\theta}) |\, \d\theta + O(\log R) \\
\le \frac1{2\pi} \int_0^{2\pi} \log |F(Re^{{\rm i}\theta}) |\, \d\theta + O(\log R) \qquad \qquad \qquad \\
= N_{\Lambda (F)}(R) + O(\log R) \le N_\Lambda (R) + O(\log R)\,, \qquad R\to\infty\,.
\end{multline*}
It remains to note that, by our assumption~\eqref{eq-N}, for each
$A<\infty$, we have
\[
\lim_{R\to\infty} \bigl (N_\Lambda (R) - 2c^{-1} R + A\log R \bigr) = -\infty\,.
\]
We arrive at a contradiction that proves the lemma. \hfill $\Box$

\end{document}